\RequirePackage{fix-cm}
\documentclass[smallextended]{svjour3}       
\smartqed  
\usepackage{graphicx}
\usepackage{epstopdf}

\usepackage{amsmath}
\usepackage{mathrsfs}
\usepackage{algorithm}
\usepackage{algorithmic}

\begin{document}

\title{Modelling power markets with multi-stage stochastic Nash equilibria }
\subtitle{ }


\author{Joaquim Dias Garcia         \and
        Raphael Chabar 
}


\institute{Joaquim Dias Garcia and Raphael Chabar\at
              Praia de Botafogo, 228 - Botafogo, Rio de Janeiro - RJ, 22250-040, Brazil \\
              Tel.: +55-21-3906-2100\\
              \email{\{joaquimgarcia\}\{chabar\}@psr-inc.com}           
                                          \\
              \\
              \hfill
              This work was presented in ICSP 2016 - Buzios, Brazil
}

\date{July 19, 2016}

\maketitle

\begin{abstract}
The modelling of modern power markets requires the representation of the following main features: (i) a stochastic dynamic decision process, with uncertainties related to renewable production and fuel costs, among others; and (ii) a game-theoretic framework that represents the strategic behaviour of multiple agents, for example in daily price bids.
These features can be in theory represented as a stochastic dynamic programming recursion, where we have a Nash equilibrium for multiple agents. However, the resulting problem is very challenging to solve.
This work presents an iterative process to solve the above problem for realistic power systems. The proposed algorithm is consist of a fixed point algorithm, in which, each step is solved via stochastic dual dynamic programming method.
The application of the proposed algorithm are illustrated in case studies with the Panama systems.
\keywords{Power Markets\and Stochastic Optimization\and SDDP\and Game-Theory\and Strategic Bidding}
\end{abstract}

\section{Introduction}\label{Intro}

Since the 80s many modifications took place in the power sectors all over the globe\cite{hu2007using}. One of the main paradigm changes was the appearance of deregulated markets which created new challenges for both industrial and academic sectors. The main modifications towards a competitive market was to replace centralized operation by energy markets in which agents are able to freely take their own decisions on both investment and production\cite{hunt2002making}. Specifics of institutional rules vary from one country to the other, however a common auction based framework is present in most of these deregulated power systems\cite{al2006electricity}\cite{maurer2011electricity}.

The first two main classes of problems arising from these changes are the investment\cite{zhu1997review} and the production problem\cite{barroso2006nash} also known as operantion and dispatch problem. While the first is deeply connected to system expansion planning and is a long term problem, with a time scale of around 10 years, the second one is a shorter term problem, of scales varying from days to few years. These problems are clearly distinct and this work is concerned with the second one, the production problem, in the specific case of hydro thermal systems.

The production problem is concerned with the agents' energy offer, or production, in a competitive market. Such problem is connected to the operation problem since the system operation will be based in those offers. In hydro thermal power systems, the centralized operation problem is already a non-trivial problem, whose solution has been well established \cite{pereira1991multi}. In many cases, its possible to reduce the first problem to the latter, also the operation problem can be used to give proxies for the solution of the production problem.

Depending of the nature of the system, hydro thermal in our case, the production problem also have some special cases which have significantly different solutions. The hydro-thermal energy production problem is considerably more complex because in inherits the basics of purely thermal problem and adds other difficulties of its own such as time coupling and stochastic processes such as incoming water to reservoirs by means of rain\cite{steeger2014optimal}. 

In a deregulated environment we have a second key difficulty in the market modelling: effect of market power \cite{hobbs2000strategic}. It may occur  that some agent, or a small set of agents, hold a significant proportion of the system's energy capacity and then be able to manipulate prices at its own will, monopolies are the extreme situation when some agents do control all the system capacity. Such situation is terrible for consumers that might be exposed to unreasonable energy prices.

\subsection{Main Contributions}\label{Contrib}

According to \cite{steeger2014optimal} the purely thermal version of the production problem has been vastly studied and a lot of good methodologies were proposed. Also, for short term hydro production problem there are a few works that devise interesting and fairly complete solutions. However, for medium term hydro production problem the literature is very limited, due to the high complexity of a problem including time coupling and market power.

Our goal is to detail a complete methodology to simulate the behaviour of a Hydro-Thermal Market. The result of the simulation procedure should approximate the production problem solution. Considering multiple agents capable of exercising market power is a common framework in game theory, whilst taking into account time coupling and uncertainties is the goal of stochastic dual dynamic programming, SDDP. Therefore we will proposed an algorithm based on the combination of SDDP and non-cooperative game theory following the lines of \cite{flach2010long}. 

However, we differ from \cite{flach2010long} because we consider price and quantity bids instead of simple quantity bids. Also the former work only considered one price maker agent and relied in the existence of price takers, while we shall not assume the existence of price takers and we will consider multiple price maker agents.

The use of algorithmic game theory will assume rationality of the agents which is a typical assumption in auditing markets for market power abuses. Thus, the method is of great use for market regulators, besides the straightforward applications in the agent energy production problem.
	
\subsection{Organization}\label{Organ}

This work is organized as follows: chapter two is concerned with the fundamental characteristics of power markets. In chapter three a few models concerning specifics of power markets well be presented. In Chapter four, the simulation algorithm is described together with the corresponding concepts of game theory. Chapter five contains the case study and a few representative results. Finally, in chapter six, conclusions are drawn together with future study topics.

\section{Power Markets}\label{market}

Power Markets have received a lot of attention in the last years, its an important economy area and moves huge amounts of money. Many countries' power sectors have been to through deep modifications creating the most varied problems for researches and practitioners. System deregulation took place in many countries creating an interesting framework for energy trade.

In those deregulated systems, competition takes place in a auction based framework  known as day ahead market \cite{krishna2009auction}\cite{maurer2011electricity}. Agents actually rely on bidding strategies in daily basis to improve their revenue. Consequently, devising methodologies to obtain good bidding strategies has become an important problem studied worldwide.

On the other hand, system operators also rely on bidding models to minimally regulate markets and avoid abusive prices in order to protect consumers\cite{borenstein2000diagnosing}\cite{wolak2014effective}. Presenting a complete methodology to simulate agents behaviour in such markets is crucial to detect market power and abuses. Simulation of energy Markets have been studied in \cite{villar2003hydrothermal} and \cite{kelman2001market}.

\subsection{Key concepts}\label{secKey}

In this Section we present some key concepts of power markets. Such concepts will be essential to understand the problem of market simulation and its main difficulties associated with the bidding problem. Most of the nomenclature used here follows the excellent literature review exposed in \cite{steeger2014optimal}.

\subsubsection{Cost based versus Bid based power systems}\label{secCostOffer}

The optimal operation problem for the case of cost based power systems is a classical and well studied problem since the famous solution by the Stochastic Dual Dyanamic Programming, SDDP, \cite{pereira1991multi} for the case of hydro thermal systems. For completeness of this work that borrows heavily from the core SDDP idea we reference the cost dispatch in the appendix 1.

In cost based these systems all the agents have operation costs and the dispatch problem minimizes the overall system operation cost without absolutely no interference of the generators. This is a centralized operation typical of regulated power systems such as Brazil and most of South and central America.

On the other hand we have the bid based model in which generators would not need to disclose their operation costs. In this model, all the generators hand their energy offer to the operator. These offers, also known as bids, are usually represented by price and quantity curves reflecting amounts of energy some generator is willing to sell for some given price.  Many countries such as United states and New Zealand have adopted such deregulated system.

\subsubsection{Day-Ahead markets}\label{secDay}

Existing offer based power markets are built up from the fundamental concept of day ahead markets\cite{rebennack2010handbook}. Every country has its own singularities on the basic rules for their markets, however, the following definition is extremely general and condensates the most important characteristics.

Following \cite{barroso2006nash}, a day ahead market consists of three phases: bidding, market clearing and pricing. The bidding phase is when the agents, relying on their own methodologies, decide their bids in form of price and quantity curves. 

Market Clearing is the process of deciding which agents will be dispatched, it is performed by the system operator and fundamentally decides the optimal dispatch by minimizing the overall cost of meeting the demand. If no extra rules, such as nodal pricing\cite{taylor2015convex}, are included in the clearing, the operator simple dispatches the agents with cheapest bids. The price of the most expensive agent dispatched is the system marginal cost, also known as spot price. This step is performed in the exact same way of the single stage thermal dispatch of appendix \ref{thermalop}, but instead of having plant's cost and capacities we have agents prices and quantities. 

Given the a demand $d$ and set of prices $p_i$ and quantities $q_i$ indexed by the set $A$ we have the simplest market clearing problem given by:

\begin{align}
& {\text{minimize}} && \sum_{i \in A}{p_i {e}_{i}} \\
& \text{subject to} \quad \notag \\
&&& \sum_{i \in A}{e_{i}}  = d_{}  && \leftarrow {\pi_p}_{} && \label{marketloadbal}\\
&&& e_{i} \leq \overline{q}_i && &&, \forall i \in A 
\end{align}

 Clearly this problem is solved simply by the merit order and yield to the system spot price $\pi$, the dual variable of the load balance constraint (\ref{marketloadbal}), see appendix 1. 

Finally, pricing is the phase in which the agents are paid for the energy they sold. All the energy is bought at spot price value no matter the position of the agent or plant in the merit order, it only matters if some of its quantity was allocated to meet demand.

The above described procedure is a simplified version and for instance transmission constraints \cite{taylor2015convex} can complicate the problem and dispatch agents without following merit order. 

\subsubsection{Price makers and price takers}\label{secPMPT}

In the cost based power systems all the agents are the same to the eye of the system operator since they simply hand in their operation cost, which is a physical characteristic. However, in offer (or bid) based markets we must depict a clear distinction between two kinds of agents. 

The first group is named \textit{price takers} because their energy offer will not affect the system spot price, because they are typically small generators compared to the overall system demand. The problem of obtaining optimal price and quantity offers for these agents is named the \textit{bidding problem}, which is very important on its own.

The second group, the \textit{price makers}, is composed of agents whose offer can indeed modify the spot price. These agents usually have generation capacity that can be compared to the system demand, or better, represent a significant proportion of the overall system capacity. For these agents the problem of obtaining optimal bids is known as \textit{strategic bidding problem}.

\subsection{Bidding Strategies}\label{secBid}

In the aforementioned excellent review of \cite{steeger2014optimal} the hydro bidding problem main characteristic are exposed so that methodologies can be separated between niches. The bidding problem will be tackled in significantly different manners depending on problem time scale, agents nature (price makers and price takers) and power plants nature (Thermals or Hydros).

Note that the Thermal bidding problem is indeed relevant for Hydro bidding because frequently hydros can be approximated locally by thermals , with respect to time. However, instead of operation variable costs, hydro plants will have opportunity cost see appendix \ref{thermalop} and \cite{pereira1985stochastic}.

In hydro modelling we can have varied time scales, we highlight 4 basic ones. \textit{Immediate Term}, represented by a minute scale where intra day operations is carried, dispatch is performed to comply with demand variations and point-wise problems, reserves and unit commitment are fundamental here. \textit{Short Term} with  a time scales ranging from days to months. \textit{Medium Term} representing one to five years range. Finally \textit{Long Term} horizon ranging from 10 to 30 years, which is a typical framework for generation expansion. Note that in the first two scales inflows are usually treated as deterministic data, whereas in the last two time horizons inflows are stochastic data.

For hydro bidding, the most relevant scales are short and medium terms. In short term, the small variability of inflows makes the problem similar to the thermal bidding problems. However, in medium term, the inflows may vary significantly as exposed in \cite{yakowitz1982dynamic}\cite{kauppi2008empirical}\cite{gjelsvik2010long}, leading to more complex problems with fewer works, mainly for the price maker case\cite{barroso2006nash}.

This reduced literature is the main motivation of this work to present a complete simulation methodology for power markets with hydro penetration. More specifically, we are focusing in the case of medium term problems with hydro price makers where literature is even more scarce.

Aiming at the hydro bidding problem with price makers we shall present the state of the art literature for the other cases in the following sections. We will present four main cases: 1) thermal price taker, 2) thermal price maker, 3) hydro price taker and our goal 4) hydro price maker. 1) and 2) are facilitated because they present no time coupling which prevents straightforward application of their results to 3) and 4). On the other hand in 1) and 3) no agent bid will affect the price leading to simpler structure and interesting solutions that do not apply directly to 2) and 4).

\subsubsection{Thermal Price Taker}\label{secTPT}

This is the simplest version of the problem presented here. We have no time coupling and no market power being exercised by agents. Due to this characteristics a fairly complete solution solution was devised by  Gross \textit{et al.} \cite{gross2000generation} by resorting on economic ideas. 

It was proved that the optimal bidding strategy is to bid the variable cost operation. Therefore the spot price resulting from the market clearing in competitive markets will converge to the prices in centralized dispatch, see appendix \ref{thermalop}, provided we do not have time coupling or price takers. This result is widely used to make price forecasting in power markets\cite{bunn2000forecasting}.

For more complete reviews on the thermal price taker version of the bidding problem the reader is referred to \cite{kwon2012optimization}\cite{conejo2001mathematical}\cite{ventosa2005electricity}.

\subsubsection{Thermal Price Maker}\label{sectTPM}

Now the problem is complicated by the presence of price makers existence. Still no time coupling is present as in last section. Price makers have the capability of affecting the spot price by modifying their bids, therefore we have a circular problem in which bids affect prices and prices affect bids. 

This circularity led to the common framework of bi-level mathematical programs. The idea is to have a first level model to optimise the agents bids and second level to perform market clearing and decide the spot prices. Fortunately the second level is a linear program (LP), thus the use of the so called Karush-Kuhn-Tucker optimality conditions\cite{boyd2004convex} to convert the minimization problem into a set of non-linear equations where the spot price is indeed a variable. By adding these equations to the first level of the optimizations problem we have a non-linear (indeed non-convex) optimization program with a special structure known as \textit{Mathematical Program with Equilbrium Constraints}, MPEC\cite{luo1996mathematical}. A lot of research have been done around this problem and the reader is referred to \cite{pereira2005strategic}\cite{hobbs2000strategic}\cite{hobbs2004modeling} for further information.

Many algorithms were proposed to solve the case of a single price maker agent, some of these are: mixed integer linear programs, MILP, \cite{pereira2005strategic}; LCP\cite{hobbs2004modeling}, heuristics \cite{conejo2002optimal}, simulations \cite{villar2003hydrothermal} and tailored procedures \cite{hobbs2000strategic}. For the multi agent case we have \cite{de2003stackelberg} where Stakelberg equilibria is achieved by a tailored algorithm, and \cite{barroso2006nash} where the MILP model of \cite{pereira2005strategic} is expanded to accommodate the multi-agent case using Nash equilibrium, see appendix 2. 

Game theoretic approaches have two interesting characteristics that will influence their applications. Firstly, in this approach the number of states can grow rapidly, which can make many practical size instances intractable \cite{cau2002co}. Secondly, players rationality is assumed, which might not be a good assumption for some applications such as bidding strategies devised by companies, although this assumption is a reasonable and common assumption on the side of the system operator. One alternative to the game theory is the scenario approach, in which a price maker company considers its rivals bid in scenario form to obtain a bidding strategy in extreme situations \cite{baillo2004optimal}.

\subsubsection{Hydro Price Taker}\label{sectHPT}

Now we relax once more the the existence of price makers and go back to systems comprising only price takers, however, now we add hydro plants to the system. Now the problem is modified once more because time coupling and stochastic inflow processes render most of the solutions of last section useless.

A seminal reference to the problem solution is given in Gjelsvik \textit{et al.} \cite{gjelsvik1999algorithm}, which addresses to the hydro price taker problem in medium term horizon. The authors consider stochastic data from both inflows and spot prices and devise a mixed SDP/SDDP algorithm. The proposed algorithm is then applied in \cite{fosso1999generation} to the Norwegian system.

Other methodologies were proposed to solve the problem in either short or medium term. Some of these methods applied to short term are: Lagrangian relaxation \cite{redondo1999short}, Neural Networks \cite{dieu2009improved}, Mixed Integer Non-linear Programming \cite{catalao2012optimal} and Mixed Integer Quadratic Programming \cite{pousinho2012scheduling}. Methods applied to medium term include Optimal Control \cite{mizyed1992operation}, Markov Decision Processes \cite{lohndorf2010optimal} and Sequential Stochastic LP \cite{helseth2010long}.

One interesting result of the application of a model similar to the one proposed in \cite{gjelsvik1999algorithm} is that if price taker hydros are located in different basins the market is indeed efficient and the dispatch of offer based systems converges to the one of cost based systems, this is shown in Lino \textit{et al.}\cite{lino2003bid}. Moreover \cite{lino2003bid} also shows how to fix the problem for agents in the same basin by creating the so called water wholesale market.

\subsubsection{Hydro Price Maker}\label{HPM}

Finally we arrive to the most complex of the version of bidding problems present in this work. Due to the presence of hydros we have to consider time coupling. Also the existence of price makers add non-convexities to the problem. Finally in case of medium term problems the stochastic inflow process adds up to the other difficulties of this problem.

Among the early references on the hydro price maker problem we find \cite{scott1996modelling}, which combines game theoretical ideas to dual dynamic programming in the framework of deregulated markets. The problem is modelled for 2 hydros as a Cournot duopoly and a few approximations are made to get rid of non-linearities. Considering stochastic inflows, \cite{kelman2001market} devises a Nash-Cournot model for the multi-stage problem with price makers solved via SDP. However both models are restricted to a small number of reservoirs due to the curse of dimensionality.

Other attempts to tackle this problem include: in medium term scale, more specifically around one year, \cite{baslis2011mid} by using MILP to obtain an operation for one price maker producer. On short term scale \cite{pousinho2013risk} considered uncertain rival bids and applied MILP to solve the strategic bidding problem and \cite{wehinger2013modeling} modelled the multi-agent electricity market by relying on model predictive bidding. 

 Finally, we highlight the work of Flach \textit{et al.} \cite{flach2010long} that models the hydro price maker problem with a single price maker in a stochastic framework applying SDDP to solve a medium term multi-period problem. As we have seen SDDP depends on the convexity of the problem to work, therefore the price maker revenue curve, which is know to have a sawtooth shape \cite{conejo2002optimal} \cite{barroso2006nash} is approximated by its concave hull. This allow the solution of problem instances with multiple reservoirs and time-steps.


 This concludes our review on the bidding problem so that we can proceed with the description of a few models that capture some specifics of the power markets. All of them will be pieces of the final simulation methodology.

\section{Power Markets Models}\label{chapModels}
	
	In this section we shall present four models that will be building blocks for the proposed simulation methodology. The first model is a standard SDDP formulation for centralized dispatch. Two of the models deal with optimal production strategies, one for price makers and other for price takers. The other model is basically a tool for converting between quantity bids to price and quantity bids.
	
	The models here presented will make heavy use of the classical SDP and SDDP algorithm and its applications to the optimal dispatch problem. The reader is referred to the Appendix for more information about the algorithms and the dispatch problem.
	
	We start with the notation used in all the upcoming models:

\subsection{Basic Notation}\label{secNotBasic}

The three models in the sequence apply to agents separately so the following sets, constants and variables are referring to a single agents.	
	
\subsubsection{Sets}\label{secNotSet}

$H$ is the set of hydro plants

\noindent 
$G$ is the set of thermal plants

\noindent
$R$ is the set of renewable energy plants (solar and wind power)

\noindent
$T$ is the set of states

\noindent
$S$ is the set of scenarios

\noindent
$L$ is the set of openings: Openings are the conditioned outcomes for the next stage

\noindent
$B$ is the set of blocks: Blocks are used to represent intra-stage varying load levels

\noindent
$M(i)$ set of hydro plants that spill or turbine to hydro $i$

\noindent
$P(i,t)$ set of autoregressive lags of hydro $i$ in stage $t$

\noindent
$\mathcal{M}$ is the set of cuts if some stage, understood by context

\subsubsection{Constants}\label{secNotCte}

$c_{t,j}$ is the unit cost of plant $j$ in stage $t$

\noindent
$\rho_{i}$ is the hydro production factor for hydro $i$, this value is used to convert from water turbined to energy generated

\noindent
$\varphi_{t,i}^p$ is the $p^{th}$ autoregressive coefficient of the inflow for plant $i$ at stage $t$

\noindent
$\overline{u}_i$ is the maximum turbined water by hydro $i$

\noindent
$\overline{v}_i$ is the maximum reservoir volume of hydro $i$

\noindent
$\overline{e}_i$ is the maximum energy generated by hydro $i$

\noindent
$\overline{g}_i$ is the installed capacity of thermal plant $i$

\noindent
$\hat{a}_{t,i}^s$ is the simulated inflow value for hydro $i$ in  stage $t$ and scenario $s$. If inside square brackets includes also all the previous values

\noindent	
$\hat{\xi}_{t,i}^l$ is the residual value for the AR process hydro $i$ in  stage $t$ and opening $l$

\noindent
$d_{t,b}$ is the system demand in stage $t$ and block $b$

\subsubsection{Decision Variables}\label{secNotVar}

We shall omit the subscript referring to the scenario $s$ for simplicity.

\noindent
$g_{t,b,i}$ energy generation of thermal plant $i$ in stage $t$ and block $b$

\noindent
$e_{t,b,i}$ energy generation of hydro plant $i$ in stage $t$ and block $b$

\noindent
$u_{t,i}$ turbined water by hydro $i$ in stage $t$

\noindent
$x_{t,i}$ spilled water by hydro $i$ in stage $t$

\noindent
$v_{t,i}$ reservoir level of hydro $i$ ate the beginning of stage $t$

\noindent
$\Delta_{t,b}$ overall produced energy by some agent at stage $t$ and block $b$

\noindent
${\pi_h}_{t,i}$ water marginal cost of hydro $i$ (a dual variable from water balance constraint)

\noindent
${\pi_p}_{t,i}$ system spot price in stage $t$ and block $b$ (a dual variable from load balance constraint)
 
\subsubsection{Coefficients of SDDP cuts}\label{secNotCuts}

The cutting planes generated during the SDDP algorithm have the parameters:

\noindent
$\varphi^m_{t}$ linear coefficient from the $m^{th}$ cut of stage $t$

\noindent
${\varphi_a}^m_{t,i}$ angular coefficient for inflow of hydro $i$ in the $m^{th}$ cut of stage $t$

\noindent
${\varphi_h}^m_{t,i}$ angular coefficient for reservoir level of hydro $i$ in the $m^{th}$ cut of stage $t$

\subsubsection{Agents Notation}\label{secNotAg}

In centralized models $G$ ,$H$ ,$R$ represent the set of plants of the whole system. Whenever we have a model applied to a single agent the sets $G$ ,$H$ ,$R$  represent the sets of plants for the current agent, that is, the one the problem is being solved for. If we have multiple agents at some point, we shall differentiate between their generators sets by a subscript, for instance, $H_i$ will represent the set of hydros of the $i^{th}$ agent and $H_{-i}$ will represent the set of hydros of all the other agents. The agents set is represented by $\mathcal{A}$.

\subsection{Cost Based Operation}\label{secSDDPcost}

For the sake of completeness and better understanding of the forth coming models we present the classical SDDP formulation for the centralized optimal dispatch problem. For further details the reader is referred to the Appendix.

The standard SDDP formulation is completely described by the sub-problem, which is solved in every stage and scenario. The following linear program represents the sub-problem of the centralized dispatch:

\begin{align}
& {\text{minimize}} & & \sum_{i \in G, b \in B}{c_j {g}_{t,b,i}} + \frac{1}{|L|}\sum_{l \in L}{\alpha_{t+1}({v}_{t+1,H},[{a}_{t+1,H}^l])} \label{sddpOBJ} \\
& \text{subject to} \quad & & \notag\\
&&& v_{t+1,i} = \hat{v}_{t,i}^s + \hat{a}_{t,i}^s - (u_{t,i} + x_{t,i}) + \sum_{j \in M(i)}(u_{t,j} + x_{t,j}) & & \leftarrow {\pi_h}_{t,i}& &,  \forall i \in H \label{sddpHydBal}   \\
&&& \sum_{b \in B}{e_{t,b,i}} = \rho_i u_{t,i} && &&, \forall i \in H \label{sddpHydCoef}\\ 
&&& \sum_{i \in H}{e_{t,b,i}} + \sum_{i \in G}{g_{t,b,i}}  = d_{t,b} - \sum_{i \in R}{r_{t,b,i}}  && \leftarrow {\pi_p}_{t,b} &&, \forall b \in B\label{sddpLoadBal} \\
&&& v_{t+1,i} \leq \overline{v}_i && &&, \forall i \in H \label{sddpUB1}\\
&&& u_{t,i} \leq \overline{u}_i && &&, \forall i \in H \\
&&& e_{t,b,i} \leq \overline{e}_i && &&, \forall i \in H \\ 
&&& g_{t,b,i} \leq \overline{g}_i && &&, \forall i \in G \label{sddpUB2}\\
&&& a^l_{t+1,i} = \sum_{p \in P(i,t)} {\phi_{t,i}^p \hat{a}_{t+1-p,i}^s }+ \hat{\xi}_{t,i}^l  && &&, \forall i \in H, l \in L \label{sddpAR}
\end{align}	

Of crucial importance in the SDDP algorithm we have the representation of the future cost function in terms of hyperplanes, which are actually included directly in the sub-problem:

\begin{align}
& \alpha_{t+1}({v}_{t,H},[{a}_{t+1,H}^l]) = \hspace{0.1in}\notag \\
& {\text{minimize}} \ \ \ \  \alpha_{t+1} \label{fcfOBJ}\\
& \text{subject to} \quad  \alpha_{t+1} \geq \varphi^{m}_{t+1} + \sum_{i \in H}{{\varphi_h}_{t+1,i}^m  {v}_{t,i}}+ \sum_{i \in H, p \in P(i,t)}{{\varphi_a}_{t+1,1,i}^m {a}_{t+1-p+1,i}^l} &&, \forall m \in \mathcal{M}, l \in L \label{fcfHyperPlanes}
\end{align}

In such model, we have the following constraints in the sub-problem: (\ref{sddpHydBal}) represents the hydro balance; (\ref{sddpHydCoef}) enforces that the energy produced by each hydro in all blocks must match the turbined water; (\ref{sddpLoadBal}) represents the load balance that enforces how much energy must be produced to meet the demand; (\ref{sddpUB1})-(\ref{sddpUB2}) are simple Upper bounds on positive variables; and finally (\ref{sddpAR}) describes the autoregressive process of inflows. 

The objective function (\ref{sddpOBJ}) is composed of two terms: the \textit{Immediate Cost Function} given by the thermal production costs and the \textit{Future Cost Function} that represents the water values since its a function of inflows and storage in stage $t+1$. The function $\alpha_{t+1}({v}_{t+1,H},[{a}_{t+1,H}^l])$ represented by (\ref{fcfOBJ})-(\ref{fcfHyperPlanes}).
	
\subsection{Price taker offer: MaxRev}\label{secMaxRevIntro}
	
As introduced before, the price taker taker optimal energy offer problem is on its own an extremely interesting problem. Moreover, it has not only been studied a lot in recent years, but also modelling was applied in industry successfully.	

In this work, we focus in the hydro-thermal power systems in medium time scale as it was defined in \cite{barroso2006nash}. Because agents can move water, thus energy, from one stage to the other via reservoir operation the time coupling in this problem is crucial. As it will be fundamental for the derivation of the market simulation model we present a methodology to solve such problem based on the SDDP algorithm, this method is fundamentally based in \cite{gjelsvik1999algorithm}.

The methodology to be presented relies deeply on the fundamental price taker characteristic: the agent has no market power or influence on the spot prices. For thermals this problem reduces to offering the operation cost of each plant. The objective will be to maximize the agent expected revenues given a set of a distribution of spot prices.
		

The solution method for maximizing a price taker expected revenue given a temporal spot price distribution will be called \textit{MaxRev}. 

The spot prices are considered as distributions over time they indeed have a stochastic processes status just like the inflow process, therefore must be modelled as a state in a SDDP type formulation. They do have some temporal structure and could be modelled as AR processes.

However, as it was shown in \cite{gjelsvik1999algorithm}, the problem turns out to be saddle shaped: convex in the reservoir and inflows states, but concave in the spot price state. The solution was to model the spot price process as a Markov process, leading to a combination of  SDDP and SDP methods.

\subsubsection{Markov Chains}\label{secMaxRevMC}

The Markov Chain scheme is very flexible and can also be used to represent  demand load growth , uncertainties in coefficients of the objective function, such as the operating costs (resulting from stochastic fuel costs).

The Markov chain is defined by a set $K$ of states, each corresponding to a cluster of spot prices,  transition probability from state \(j\) in stage \(t\) to state \(m\) in stage \(t + 1\) is represented by \(p_{t}^{jm}\).

Spot prices are obtained from some dispatch methodology and, consequently, correspond to some scenario $s$ and stage $t$ matching, in particular, some inflow vector $a_{t,s,H}$. These spot prices are clustered in order to create the states of the Markov Chain so that each state $k$ corresponds to a set ${M}_{t}^{k} $ of spot price scenarios $\{\hat{\pi}_{t,s,b} | s \in S^k \}$, where $S^{k,t}$ is the set of scenarios that are in cluster $k$ at stage $t$. The clustering method can be of user choice although some methods such as k-means are almost the standard. The transition probabilities $p_t^km$ are estimated by counting how many simulation scenarios $s$ which are in state (cluster) $k$ in stage $t$ belong to state (cluster) $m$ in stage $t+1$.

Each cluster \(k \in K \) in each stage $t$ contains \({M}_{t}^{k}\) values of spot prices. The state variable will correspond to a cluster instead of a actual spot price value. Therefore, instead of a future cost function $\beta_t(\hat{v}_{t,H}^s,[\hat{a}_{t,H}^s],\pi^k_{t,B})$ we actually have $\beta_t(\hat{v}_{t,H}^s,[\hat{a}_{t,H}^s],k(s)) = \beta_t^{k(s)}(\hat{v}_{t,H}^s,[\hat{a}_{t,H}^s]) $, where with some abuse of notation $k(s)$ is a function that maps scenarios to clusters.

\subsubsection{Problem Formulation}\label{secMaxRevForm}

Now we present the basic model formulation for the MaxRev. Following a classic dynamic programming modelling style we go straight to the decomposed version of the problem. We have the following sub-problem (or bellman recursion):

\begin{align}
&\beta_t^{k(s)}(\hat{v}_{t,H}^s,[\hat{a}_{t,H}^s]) = \\
& {\text{minimize}}  \ \ \ -\sum_{b \in B}\pi^s_{t,b}\Delta_{t,b}+\sum_{i \in G, b \in B}{c_j {g}_{t,b,i}} + \frac{1}{|L|}\sum_{l \in L}{\beta_{t+1}^{k(l)}({v}_{t+1,H},[{a}_{t+1,H}^l])} \\
& \text{subject to} \quad & &\notag\\
& v_{t+1,i} = \hat{v}_{t,i}^s + \hat{a}_{t,i}^s - (u_{t,i} + x_{t,i}) + \sum_{j \in M(i)}(u_{t,j} + x_{t,j}) & & & &,  \forall i \in H \label{mr_hydbal}    \\
& \sum_{b \in B}{e_{t,b,i}} = \rho_i u_{t,i} && &&, \forall i \in H \\ 
& \sum_{i \in H}{e_{t,b,i}} + \sum_{i \in G}{g_{t,b,i}} +\sum_{i \in R}{r_{t,b,i}} = \Delta_{t,b}    && &&, \forall b \in B \label{mr_loadbal} \\
& v_{t+1,i} \leq \overline{v}_i && &&, \forall i \in H \\
& u_{t,i} \leq \overline{u}_i && &&, \forall i \in H \\
& e_{t,b,i} \leq \overline{e}_i && &&, \forall i \in H \\ 
& g_{t,b,i} \leq \overline{g}_i && &&, \forall i \in G \\
& a^l_{t+1,i} = \sum_{p \in P(i,t)} {\phi_{t,i}^p \hat{a}_{t+1-p,i}^s }+ \hat{\xi}_{t,i}^l  && &&, \forall i \in H, l \in L \label{mr_AR}
\end{align}

The first observation is that we have the new decision variable $\Delta_{t,b}$ that represents the amount of energy sold by the agent at the spot price of the current cluster: $\pi_{t,b}^s$. The second basic observation is that the load balance equation of the model in section (\ref{sddpLoadBal}) was replaced by equation (\ref{mr_loadbal}) which states that the $\Delta_{t,b}$ is indeed the total generation of the agent.

Observe that the spot price here is the one matching the inflow scenarios, it actually is not a necessary condition because any spot price from the correct cluster can be used, provided we choose it randomly.

Now one can simply apply the SDDP algorithm (see appendix) with the following caveats: since now we have spot price clusters as states the benders cuts generated by some sub-problem get the label of the cluster from which it was generated. This label is used to associate the cut with the correct future cost function. In scenario $s$ we are in cluster $k(s)$ and therefore we may transit to other cluster according to the process transition matrix. These possible transitions are reflected in the proportions of the opening, that is, if we have probability $p$ of transitioning to state $k2$ than $100p$ percent of the opening should be from that cluster.

The perfect match between proportions and opening numbers is hard to obtain therefore we proceed as follows: the cluster to which each opening belongs to is chosen randomly from the distribution of transition probabilities.

The cuts are defined as follows:

\begin{align}
& \beta_{t + 1}^{l,k} \geq \varphi^{m}_{t+1} + \sum_{i \in H}{{\varphi_h}_{t+1,i}^m  {v}_{t,i}}+ \sum_{i \in H, p \in P(i,t)}{{\varphi_a}_{t+1,1,i}^m {a}_{t+1-p+1,i}^l} && &&, \forall l \in L, m \in \mathcal{M}(l) \label{mr_FBF}
\end{align}

where $\mathcal{M}(l)$ is the set of cuts from the cluster that the opening $l$ belongs to.

\subsection{Detailing Bids}\label{secOptBidIntro11}

Although MaxRev can determine the optimal revenue-maximization strategy of a price-taker agent in a power market, it does not represent the specific details of the day-ahead bidding under uncertainty. We shall describe a model to perform this operation, we name it OptBid. OptBid uses the expected future revenue functions produced by MaxRev's stochastic recursion algorithm, together with auction-specific data, to calculate the price \& quantity bids of each agent in each stage (week or day) and scenario. 
		
\subsubsection{OptBid}\label{secOptBidIntro2}
		
We will decide the agent's bid in terms of price and quantity pairs as follows. Given a discrete set of prices we optimise the quantity of energy associated to each price, for each scenario, stage and block. For each inflow scenario, the bidding is carried out under uncertainty regarding the market spot prices, because for each cluster we consider a probability distribution (the spot price scenarios come from the same stochastic price model used by MaxRev in the previous run). 

The objective is to maximize the risk-adjusted sum of net revenues for that stage plus the expected revenue from the auctions in the next stages (given by the future revenue functions, also calculated by MaxRev). The use of future revenue functions ensures the adequate trade-off between the bidding revenues in a given stage and the revenues in the next stages. Otherwise, the storage devices such as hydro plants with reservoirs, pumped storage), fuel reservoirs and large-scale batteries would be emptied in the first stages. 

\subsubsection{Intuition: filling the boxes}\label{secOptBidIntuition}

Since we are associating quantities to given prices, we can think of the algorithm as an allocation of energy quantities in boxes labelled with prices. Consider the example that we have three boxes with prices 10, 20 and 30. We put energy in the box if we accept selling energy for any price above that. Suppose the current cluster spot price distribution is represented by the equiproportional values 5, 15, 25 and 35, also suppose we one thermal plant with cost 15 and capacity 10.

First of all, we put any amount of energy in box one we lose money because our cost is higher than our price, and thus we only put energy in boxes if the spot is higher that 20. If the spot is 25 we want to put some energy in the box with price 20 but not in the one with price 30, and if the spot is 35 we want to put energy in the box of price 30. But how much energy we put in which box since we are limited to 30 units? 

Since we have no revenue for spots 5 and 15, our problem is reduced to maximizing $1/4((0)+(0)+ ((10-15)q_{10}+(20-15) q_{20})+ ((10-15)q_{10}+(20-15) q_{20}+(30-15) q_{30}))$ subject to $q_{10}+q_{20}+ q_{30} \leq 10$, the we would decide to put all energy in box of price 30 because it maximizes our (risk neutral) expected profit. 

The problem get more and more complex as we have many more plants, with different prices and capacities, and even more by adding hydros and their future revenue functions.

\subsubsection{Problem formulation}\label{secOptBidForm}

The optimization is carried out separately for each agent, stage and scenario of the previous optimization. For a simpler notation, we will omit in the formulation below the indices of the agents and of the parameters that define the inflow scenarios ($s$ and ${\hat{k}}^{s}$).\\

\noindent
\( b \in B \) intra-stage blocks (typically hours)

\noindent
\( n \in N\) bid segments, or "boxes", in each hour

\noindent
\( {\hat{\Pi}}_{b,n} \) (pre-defined) bid price for box \(n\), block \(b\)

\noindent
\(q_{b,n}\) bid energy amount (decision variable) of bid box \(n\), hour \(b\)

\noindent
\( k \in \mathcal{K}(s) \) spot price scenarios from cluster $K(s)$ of current scenario $s$ in stage \(t\) (as seen in the previous models, these scenarios are obtained from the clustering of short-run marginal costs in the construction of the spot price Markov chain). Note the notation abuse, the "function $\mathcal{K}$" maps a inflow scenario, $s$, into the spot scenarios set corresponding to the spot cluster $k(s)$, i.e. if we map it spot cluster to its set of scenarios trough $f$, we have $\mathcal{K}(s) = f(k(s))$ 

\noindent
\( \delta_{b}^{k} \) total bid energy amount (decision variable) of hour \(b\), spot price scenario \( k \).

\noindent
\( {\hat{\pi}}_{b}^{k} \) spot price of hour \(b\), scenario \(k\)

\noindent
\( {\hat{\phi}}_{b,n}^{k} \) (pre-calculated) indicator function of the acceptance/rejection of bid price \( {\hat{\Pi}}_{b,n}\) for scenario \(k\): \({\hat{\Pi}}_{b,n} \leq {\hat{\pi}}_{b}^{k} \Rightarrow {\hat{\phi}}_{b,n}^{k} = 1\); otherwise \({\hat{\phi}}_{b,n}^{k} = 0\).\\

Now we ha the following problem:

\begin{align}
& \text{Max\ } \ \ \  \frac{1}{|\mathcal{K}(s)|}\sum_{k\in \mathcal{K}(s)}^{}\left\{ \left\lbrack \sum_{b\in B}^{}{\sum_{n\in N}^{}{\left( {\hat{\Pi}}_{b,n} \times {\hat{\phi}}_{b,n}^{k} \right) \times q_{b,n}}} \right\rbrack - z^{k} + \frac{1}{|L|}\sum_{l}^{}\beta_{t + 1}^{l,k} \right\} \label{ob_obj}\\
& \text{subject to} \quad & &\notag\\
& \delta_{b}^{k} = \sum_{n\in N}^{}{\hat{\phi}}_{b,n}^{k} \times q_{b,n}\label{ob_e1} \\
& \delta_{b}^{k} = \sum_{j \in J}^{}g_{b,j}^{k} + \sum_{i \in H}^{}e_{b,i}^{k} + \sum_{r \in R}^{}{\hat{r}}_{b,r}^{k} \label{ob_e2}\\
& z^{k} = \sum_{b}^{}{\sum_{j \in G}^{}{c_{j}g_{b,j}^{k}}} \label{ob_coster}\\
& v_{t+1,i}^{k} = \hat{v}_{t,i}^s + \hat{a}_{t,i}^s - (u_{t,i}^{k} + x_{t,i}^{k}) + \sum_{j \in M(i)}(u_{t,j}^{k} + x_{t,j}^{k}) & & & &,  \forall i \in H    \\
& \sum_{b \in B}{e_{t,b,i}^{k}} = \rho_i u_{t,i}^{k} && &&, \forall i \in H \\ 
& v_{t+1,i}^{k} \leq \overline{v}_i && &&, \forall i \in H \\
& u_{t,i}^{k} \leq \overline{u}_i && &&, \forall i \in H \\
& e_{t,b,i}^{k} \leq \overline{e}_i && &&, \forall i \in H \\ 
& g_{t,b,i}^{k} \leq \overline{g}_i && &&, \forall i \in G \\
& a^l_{t+1,i} = \sum_{p \in P(i,t)} {\phi_{t,i}^p \hat{a}_{t+1-p,i}^s }+ \hat{\xi}_{t,i}^l  && &&, \forall i \in H, l \in L \label{ob_inflow}\\
& \beta_{t + 1}^{l,k} \leq \varphi^{m}_{t+1} + \sum_{i \in H}{{\varphi_h}_{t+1,i}^m  {v}_{t,i}^{k}}+ \sum_{i \in H, p \in P(i,t)}{{\varphi_a}_{t+1,1,i}^m {a}_{t+1-p+1,i}^l} && &&, \forall l \in L, m \in \mathscr{M}(l) \label{ob_FBF}
\end{align}

Clearly with the exception of equation (\ref{ob_inflow}) all the other inequalities are written $\forall k \in \mathcal{K}(s)$. The Load Balance is replaced by equations (\ref{ob_e1}) and (\ref{ob_e2}) impose that for all block $b$ and spot scenario $k$ the quantities of energy allocated to price boxes are smaller the same as the energy produced, which is bounded by the physical constraints. Equation (\ref{ob_coster}) sums the thermal cost to be considered in the objective. 

Note that we are optimizing with respect to the average value and no risk aversion is being considered, but that could be changed by modifying the objective function.

\subsection{Single agent price maker offer}\label{secNashBidIntro1}

In this section we describe the basics of the algorithm proposed in \cite{flach2010long} with a few modifications, for future reference we name it NashBid.

\subsubsection{NashBid}\label{secNashBidIntro2}

The NashBid model is used to represent the strategic behaviour of price maker agents . As the tool name implies, it is based on game-theoretic concepts, namely a stochastic multi-stage Nash equilibrium, this will be clear in the simulation algorithm section. 

In this model we assume that each agent $i$ do not have complete access to the operation and production of the remaining agents, instead it is only assumed to know the set of {[}price; quantity{]} bids of the other agents. In order to proceed, these bids are converted into a set of \emph{virtual thermal plants} where unit operating cost corresponds to the bid price and the generation capacity corresponds to the bid quantity. These virtual plants will become a ``price taker'' group in the calculation of the ``price maker'' strategy of agent $i$.

This price maker strategy is determined by a multi-stage stochastic optimization algorithm, where the objective is to maximize the product amount of energy offered by the price maker and the spot price, known as agent revenue, which depends on both the amount offered by the price maker and on the virtual thermal plants.

However as we discussed in chapter \ref{market}, the aforementioned agent revenue is sawtooth shaped curve \cite{conejo2002optimal}, which is non-convex. In order to apply SDDP we follow the approach of \cite{flach2010long} and overestimate it by its concave hull.

Just like MaxRev, NashBid determines the optimal operation strategy (but now for a price maker agent). The detailed bidding strategy can be also determined by the OptBid simply by using the pertinent future benefit functions.

\subsubsection{Agent Revenue}\label{secNashBidRev}

The revenue function of agent $i$ is constructed as follows:

We start by representing all the other agents by their price and quantity offers. Therefore, we have a set $O_{-i}$ of all offers of all the agents so that we have a market clearing with respect to a offered quantity $e$ of the agent $i$ given by the following parametric linear program:

\begin{align}
z(e) = \ \ \  & \text{Min} && \sum_{j \in O_{-i} }p_j q_{j}  \\
& \text{subject to} \quad & & \notag\\
&&& \sum_{j \in O_{-i} }q_{j} = d-e &&\leftarrow \pi_d(e) && \\ 
&&& q_{j} \leq \overline{q}_j && &&, \forall j \in O_{-i} 
\end{align}

At the optimal solution we have the spot price $\pi_d(e)$, and thus the profit of agent $i$ is given by $\hat{R}(e) = e\times \pi_d(e) $, which is sawtooth shaped, which could be represented  perfectly in a MILP. Since we are about to apply SDDP we must obtain the concave-hull $R(e)$ as in \cite{flach2010long}.

This concave hull can be represented by finite set of points $(e_i,R(e_i)), \forall i \in \mathcal{R}$, whose convex combinations form a polyhedron that represents the hypograph of $R$.

\begin{figure}[h!]
\centering
\includegraphics[width=0.7\linewidth]{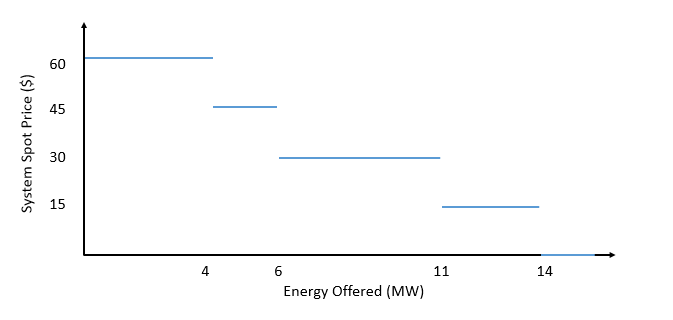}
\caption{Spot Price as a function of energy quantity offer for some price maker agent}\label{figRevFunc1}
\end{figure}

Now we exemplify the construction of this revenue function we a few simple figures. Firstly, Figure \ref{figRevFunc1} shows the system spot price as a function of the energy offered by some agent. Note that the function is piecewise constant, because the spot price is only altered when the next bid is displaced with more energy offered by the price maker agent.

\begin{figure}[h!]
\centering
\includegraphics[width=0.7\linewidth]{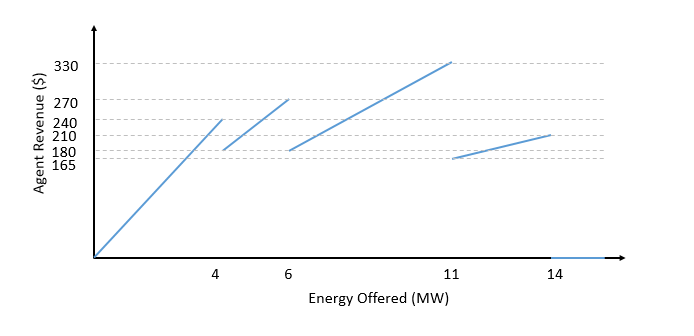}
\caption{Revenue function of a price maker agent}\label{figRevFunc2}
\end{figure}

Secondly, we have to multiply the spot price by the energy quantity offer to obtain the revenue as a function of the energy offer. The obtained sawtooth shaped function is presented in Figure \ref{figRevFunc2}.

\begin{figure}[h!]
\centering
\includegraphics[width=0.7\linewidth]{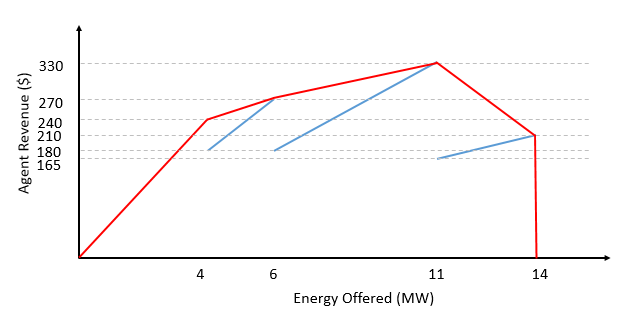}
\caption{Concave hull of the Revenue function , in red, original function in blue}\label{figRevFunc3}
\end{figure}

Finally, we obtain the concave hull of the revenue function, this is presented in Figure \ref{figRevFunc3}.

\pagebreak
	
\subsubsection{Problem formulation}\label{secNashBidForm}

After defining the revenue function the SDDP sub-problem construction is a straightforward modification from the MaxRev problem. Just one term in the objective function is modified, the profit term is replaced from a price taker profit whose energy offer does affect the spot price by a revenue function that represents the profit as function $R_{t,b}$.
	
\begin{align}
&\beta_t^{k(s)}(\hat{v}_{t,H}^s,[\hat{a}_{t,H}^s]) = \\
& {\text{minimize}}  \ \ \ -\sum_{b \in B} R_{t,b}(\Delta_{t,b})+\sum_{i \in G, b \in B}{c_j {g}_{t,b,i}} + \frac{1}{|L|}\sum_{l \in L}{\beta_{t+1}^{k(l)}({v}_{t+1,H},[{a}_{t+1,H}^l])} \\
& \text{subject to} \quad & &\notag\\
& \hspace{2in}\text{(\ref{mr_hydbal})-(\ref{mr_AR})} \notag
\end{align}

We have, just like in MaxRev, a markov chain process, but instead of a linear revenue function we have a more complex concave function. The remainder of the problem is identical.

\section{Proposed Market Equilibrium Model}\label{chapPorposed}

In this section we present 
a simulation algorithm that will make extensive use of the previously described models. Our goal is to simulate a multi-agent power market with presence of hydro plants and price maker agents in a medium term horizon. The most challenging characteristic of the problem are time-coupling and stochastic inflows inherited from the hydro operation and non-convexities from the price maker bids. Any attempt to model the multi-stage stochastic problem without a decomposition makes the problem intractable and a simple application of SDDP is useless due to non-convexities. Dealing with non-convexities in SDDP is a recent research field\cite{thome2013non}\cite{philpott2016midas}.

We shall build upon NashBid and, thus, rely on the work of Flach \textit{et al.} \cite{flach2010long} for the core algorithm. However instead of single price maker agent we will allow for multiple agents that can exercise market power, this will be done by applying basic ideas of non-cooperative game theory, in particular Nash equilibria.

\subsection{Algorithm overview}\label{secAlgOver}

The basic idea is to loop trough agents until no agent desire to change its strategy, therefore reaching a Nash Equilibrium. We will present 
a methodology of equilibrium in which
each agent's operation strategy is decided for all stages considering the bids of the other agents, then the operation of such agent is converted into price and quantity bids so that we can move to the next agent. This is repeated until no agent has incentive to change their bidding strategy.

\subsection{global equilibrium method}\label{secGlobalEq}

We start by assuming that each agent has a initial set of bids that came from a set of spot prices. We start with agent $1$ and build its revenue function, as in section \ref{secNashBidIntro1}, for all stages, blocks and scenarios from the information of the bid of the remaining agents. Then a complete NashBid SDDP recursion is run until convergence. From this problem solution we obtain sets of Future Revenue Functions that can be used by OptBid to update agent $1$ bid for all stages, scenarios and blocks. Now we go to the next agent and apply the same procedure until no agent changes its bid any more and multi-stage Nash equilibria is said to take place.

This procedure is highlighted in the next diagram:

\begin{figure}[h!]
\centering
\includegraphics[width=0.7\linewidth]{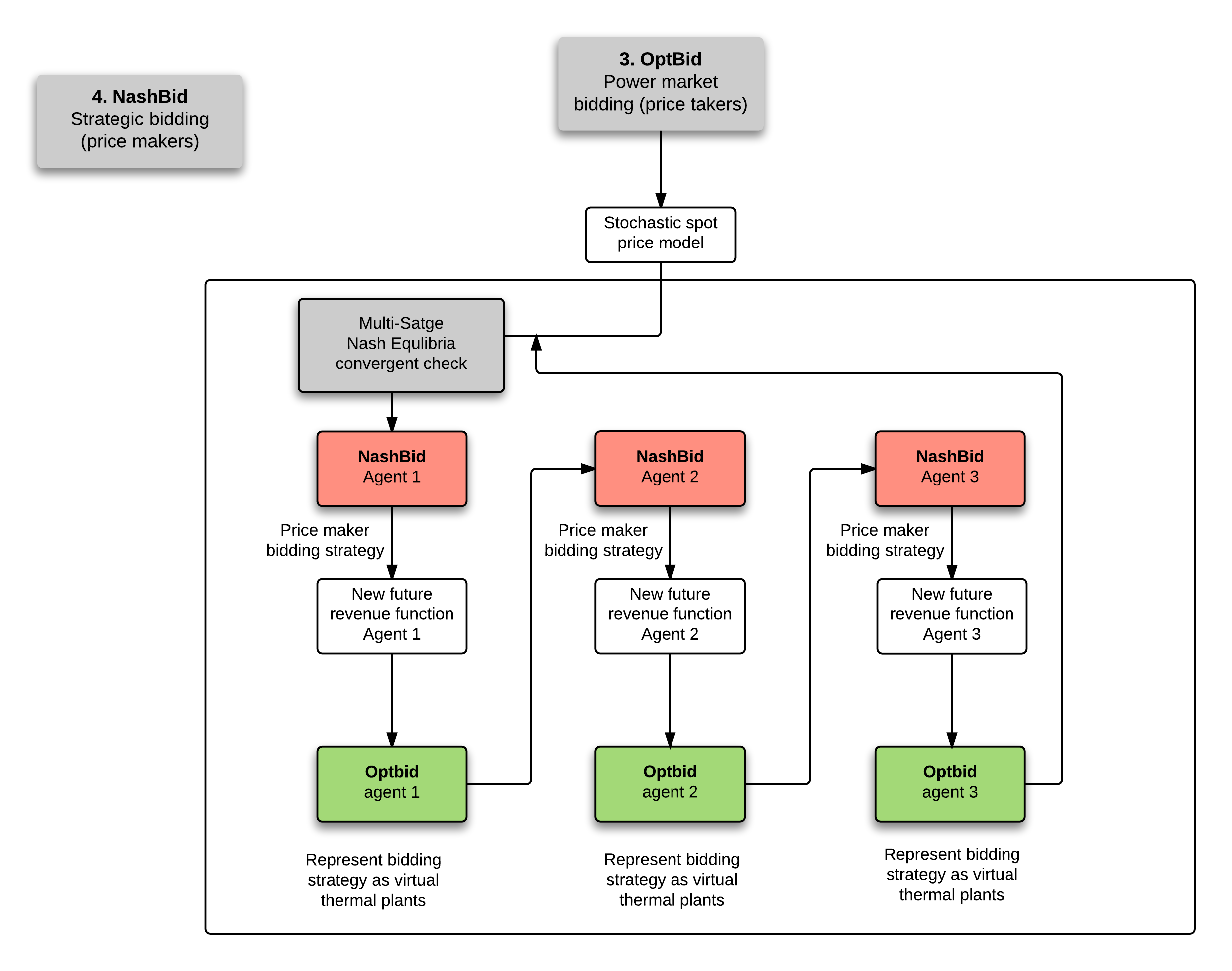}
\caption{global equilibrium method}
\end{figure}

\pagebreak

%
%
%
%

\pagebreak

\subsection{Initialization}

The methodology described to achieve Nash equilibrium can be significantly time consuming because the fixed point method may need a few iterations to converge. Having to converge many SDDP's for each fixed point iteration could take a long time just like having to converge Nash equilibrium in forward simulations.

In order to tackle this problem we propose a initialization methodology that only assumes as input the system characteristics:

Everything starts by optimising the centralized operation to obtain a first proxy of spot prices in the system. The second step is to convert the spot price scenarios into clusters with some algorithm such as k-means. With the clustered spot prices one can apply the MaxRev model to all agents and obtain their Future Revenue Functions. Given these FRFs and the spot prices we apply OptBid to all agents to obtain their bids and, subsequently, we apply market clearing process to obtain a second estimate of spot price scenarios. These spot price scenarios together with the OptBid results are used as input

This initialization procedure is summarized in the following figure:

\begin{figure}[h!]
\centering
\includegraphics[width=0.7\linewidth]{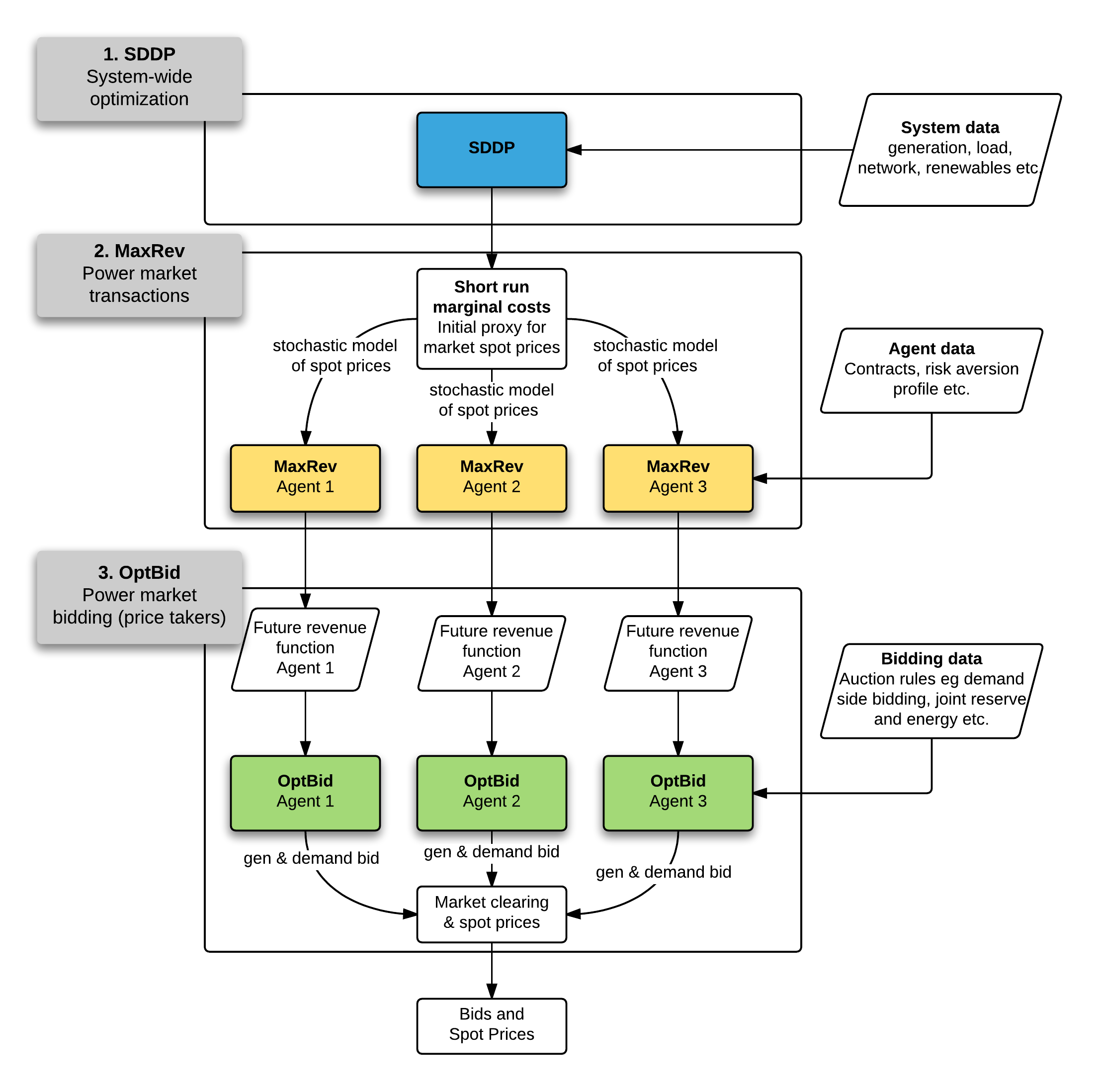}
\caption{initialization}
\end{figure}

\pagebreak

\section{Case Study}

In this section we outline basic results of simulation by applying the methodology in the system of Panama.

\subsection{The Panama system}

The Panama power system is  an hydro thermal system therefore it  includes the core difficulties that we seek to model. The system configuration is slightly different from the real system to better provide insight over the results.

The system we used in the study is composed of 22 existing thermal plants with  varied fuels such as carbon, bunker and diesel, the overall thermal power is around 1145 MW. There is 42 hydro plants with total installed capacity of 1674 MW, most of the plants are run of river and three of them do have storages.

The study is carried over 4 years with monthly resolution, in which demand varies from around 850 MW per month in the first year to around 1050 MW per month in the last year.
For the sake of simplicity we consider a single block and an autoregressive model of maximum order one for the inflows stochastic process.

Finally, in order to model a competitive market with price makers we divide the set of plants in 4 sets. The first three sets are price maker agents, each with one of the water reservoirs, and the first agent owns a cascade composed of nine plants, summing a maximum installed capacity of 681MW. All the agents have similar installed capacity of 539MW and 499 MW. The fourth set contains the remaining plants, including all the renewable energy plants, thermals and few hydros with overall capacity of 1100MW, these will all be price takers.

\subsection{Result analysis}

\begin{figure}[h!]
\centering
\includegraphics[width=0.8\linewidth]{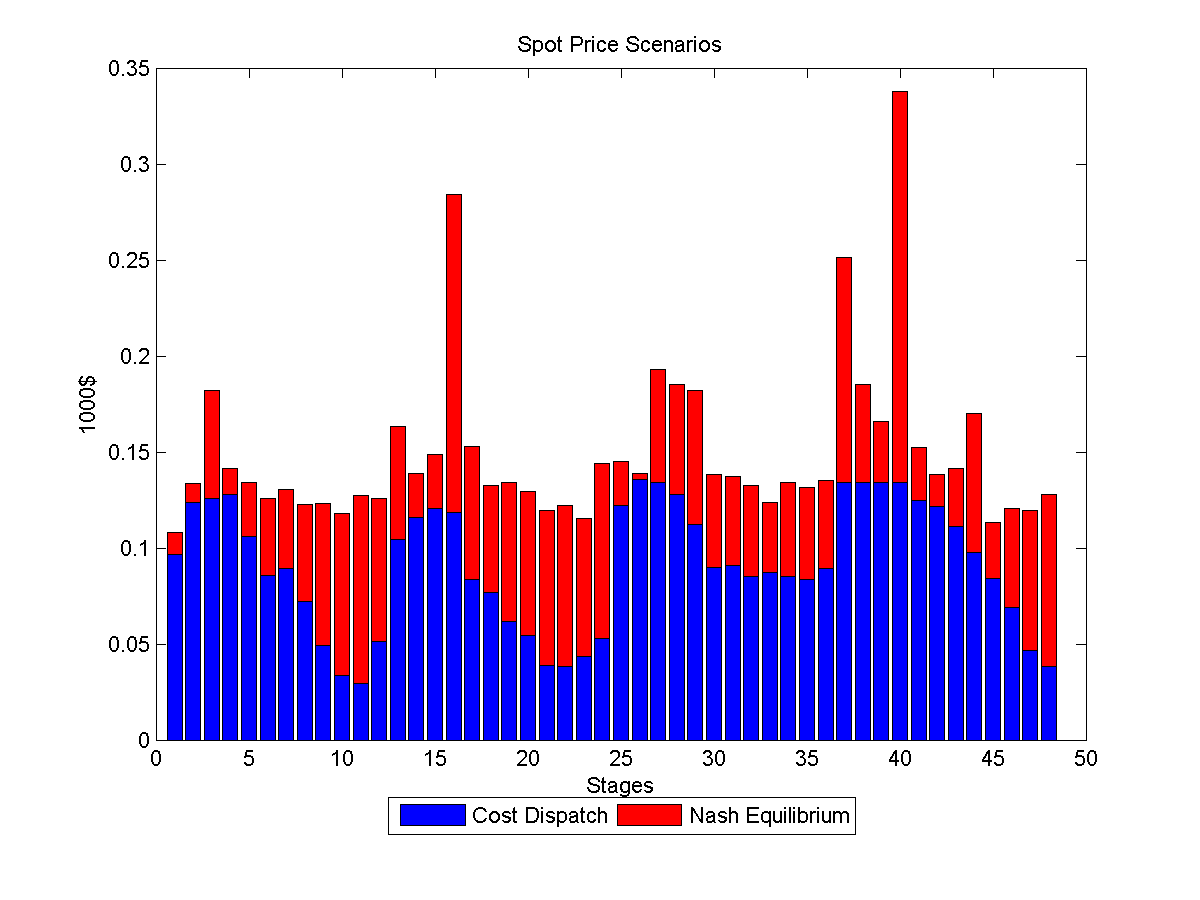}
\caption{Spot Price Mean}\label{figSpotAvg}
\end{figure}

The resulting simulation of the deregulated energy market is summarized in  the figures presented in this section. As a first illustrative figure of the result we present the monthly average of the spot price in Panama in Fig. \ref{figSpotAvg}. In blue we have the spot prices from the centralized cost minimization problem, clearly the prices are very seasonal due to the significant presence of hydro plants and the reduced number of reservoirs, only three. The average spot price in the first two years is around $\$70$, in the following to years the spot average raises due to the lack of reinforcements in the system, that is, no new plants are starting to operate after the second year. In the remainder of this section \textit{CD} will stand for Centralized Dispatch while \textit{NE} will stand for Nash Equilibrium, the results of the deregulated market simulation.

Also in Fig. \ref{figSpotAvg} we have, in red, the resulting spot prices of the deregulated system. As expected we see the prices getting significantly higher because a huge amount of the power system was allocated to price maker agents. Disregarding the spot price spikes we can still see a slight seasonality in the prices. The price spikes occur in stages where the cost dispatch attained its highest values and, thus, without too much effort the price maker agents can displace thermal plants and make this spot much higher.

\begin{figure}[h!]
\centering
\includegraphics[width=0.8\linewidth]{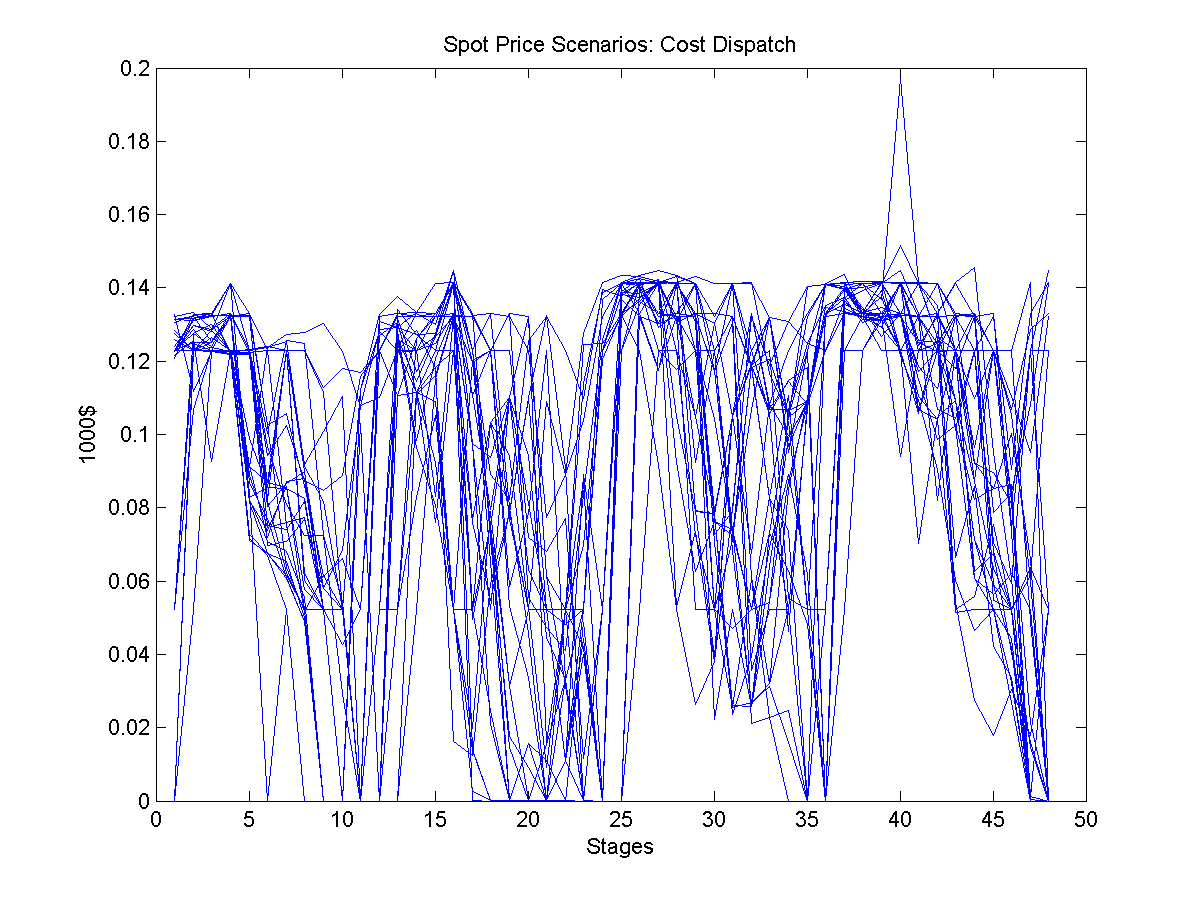}
\caption{Spot Price Scenarios}\label{figSpotScen_cd}
\end{figure}

\begin{figure}[h!]
\centering
\includegraphics[width=0.8\linewidth]{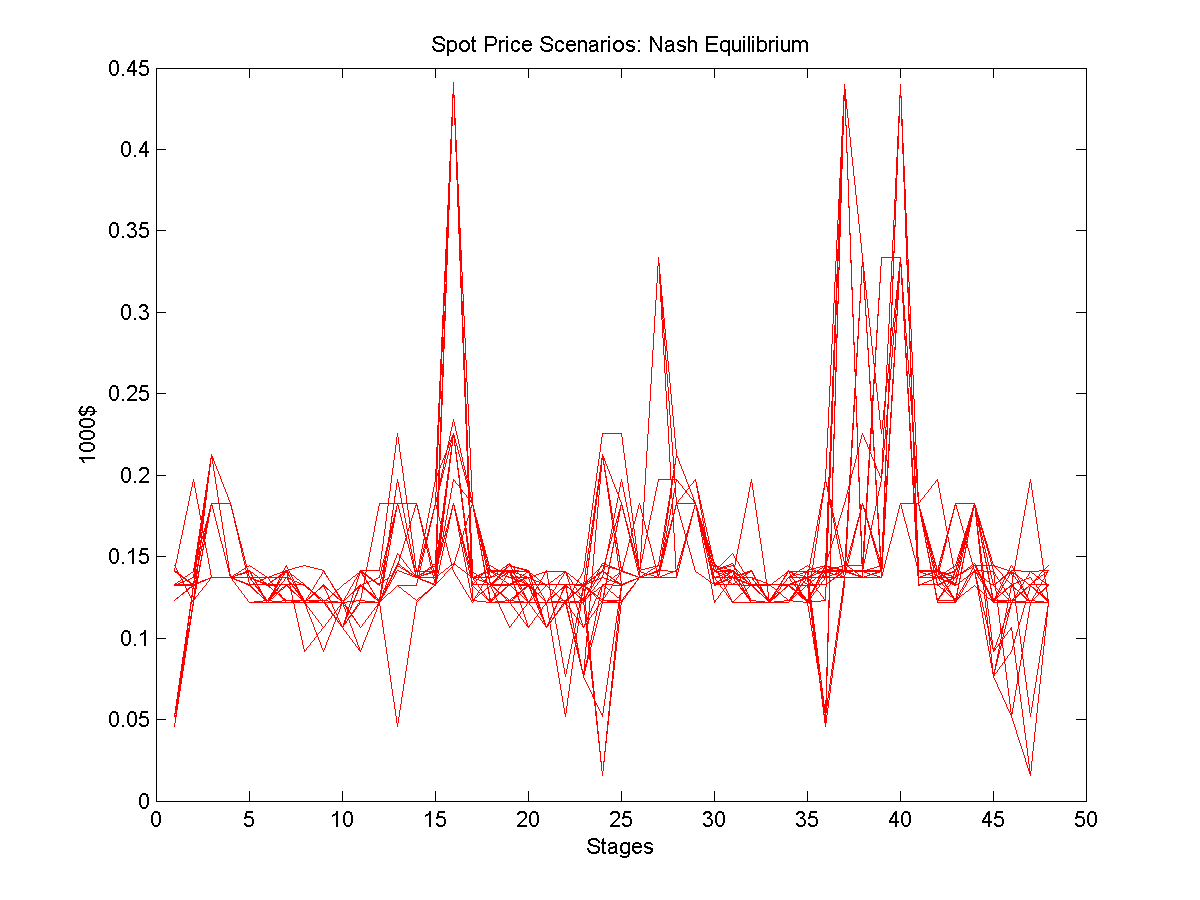}
\caption{Spot Price Scenarios}\label{figSpotScen_ne}
\end{figure}

\begin{figure}[h!]
\centering
\includegraphics[width=0.8\linewidth]{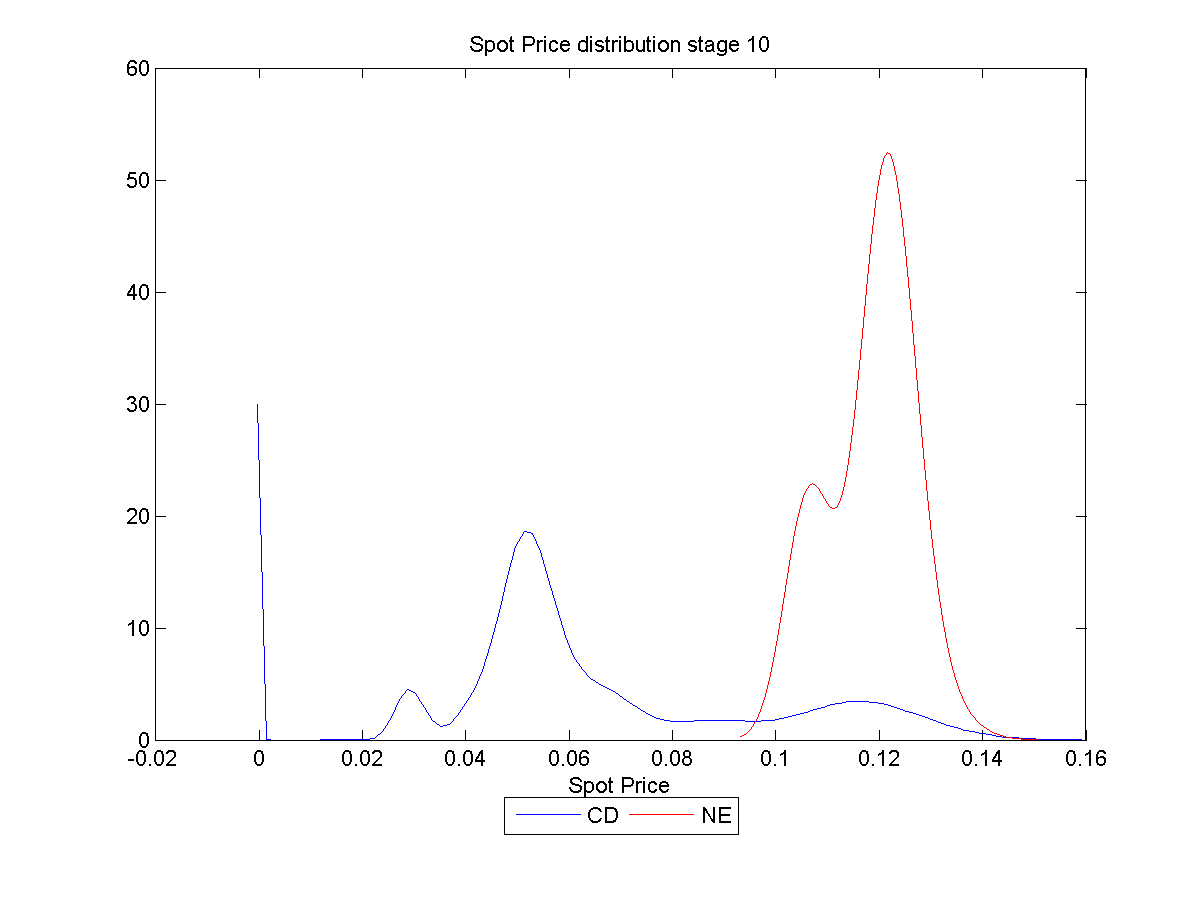}
\caption{Typical Cheap Spot Price Distribution}\label{figSpotCheap}
\end{figure}

\begin{figure}[h!]
\centering
\includegraphics[width=0.8\linewidth]{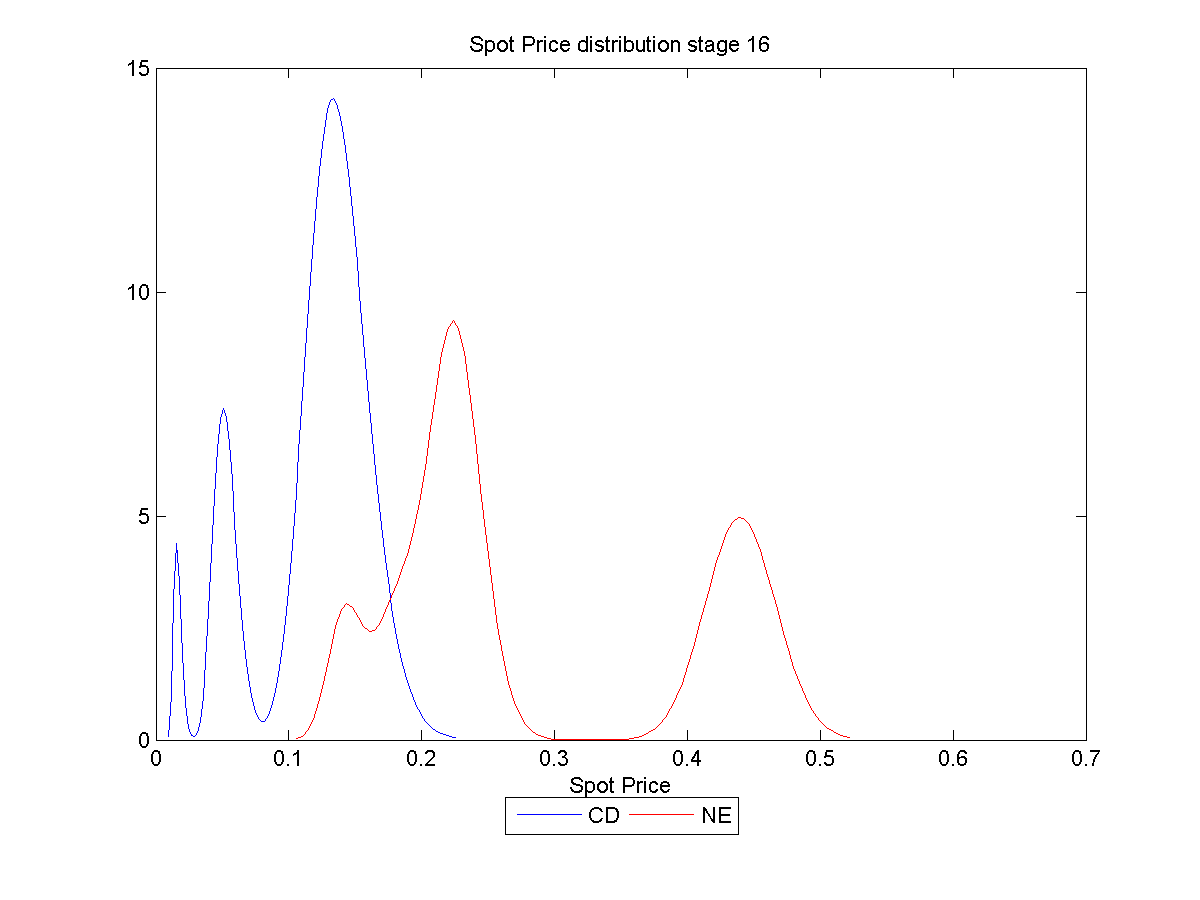}
\caption{Typical Expensive Spot Price Distribution}\label{figSpotExp}
\end{figure}

For completeness, we also present Fig. \ref{figSpotScen_cd} and Fig. \ref{figSpotScen_ne}  which shows the spot price scenarios for both cost dispatch and Nash equilibrium. From this graph we can have an idea of the variability in each stage. We also present Fig. \ref{figSpotCheap} and Fig. \ref{figSpotExp} which depict respectively typical cheap and expensive spot price distributions. As seen in the previous figures the spot price average is displaced to the right after deregulated equilibrium was attained. In cheap stages the spot is basically displaced to the highest spot price and no more that that. However, on expensive stages, the price is slightly displaced in many scenarios and in very critical scenarios the combination of high demand, small inflows and market power lead to extremely high spot prices.

\clearpage

%

\begin{figure}[h!]
\centering
\includegraphics[width=0.7\linewidth]{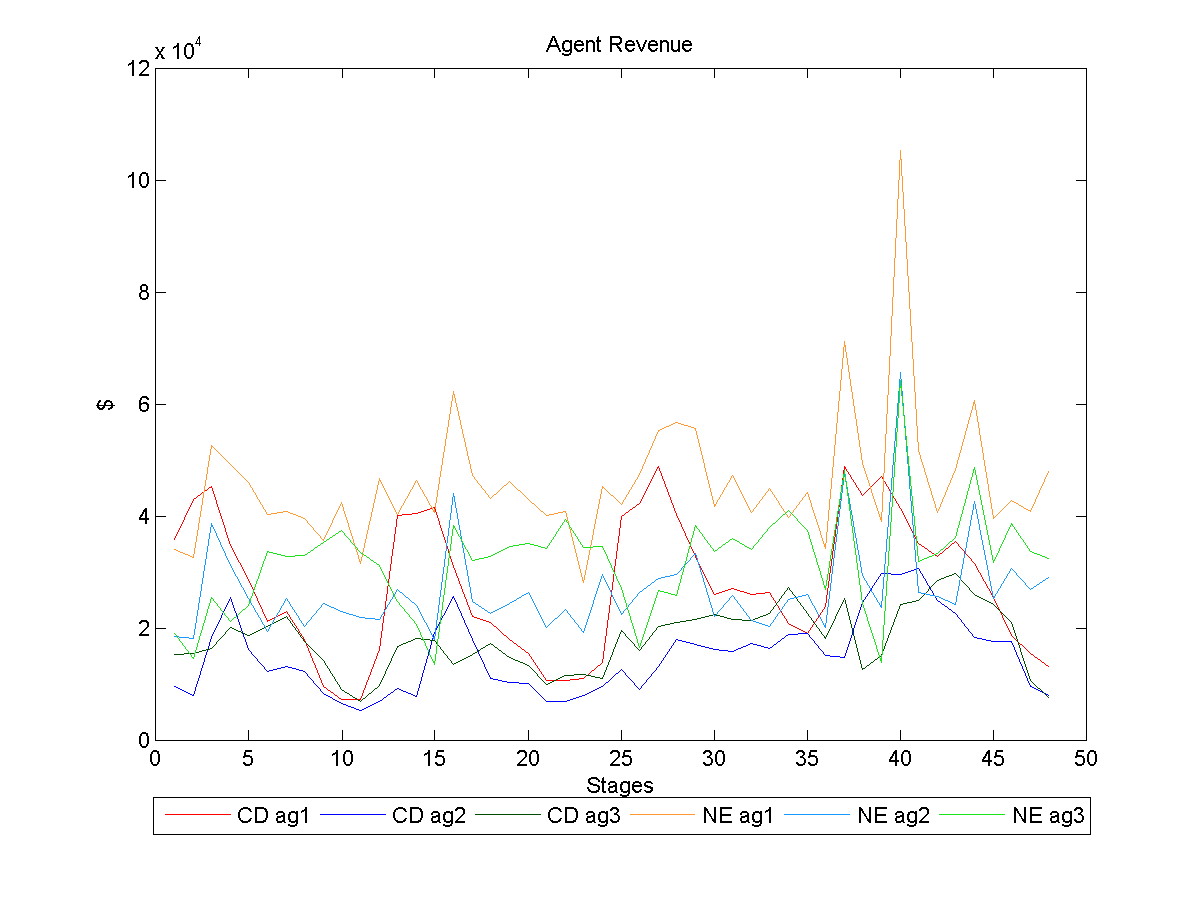}
\caption{Mean Revenue per agent}\label{figAgentRev}
\end{figure}

\begin{figure}[h!]
\centering
\includegraphics[width=0.7\linewidth]{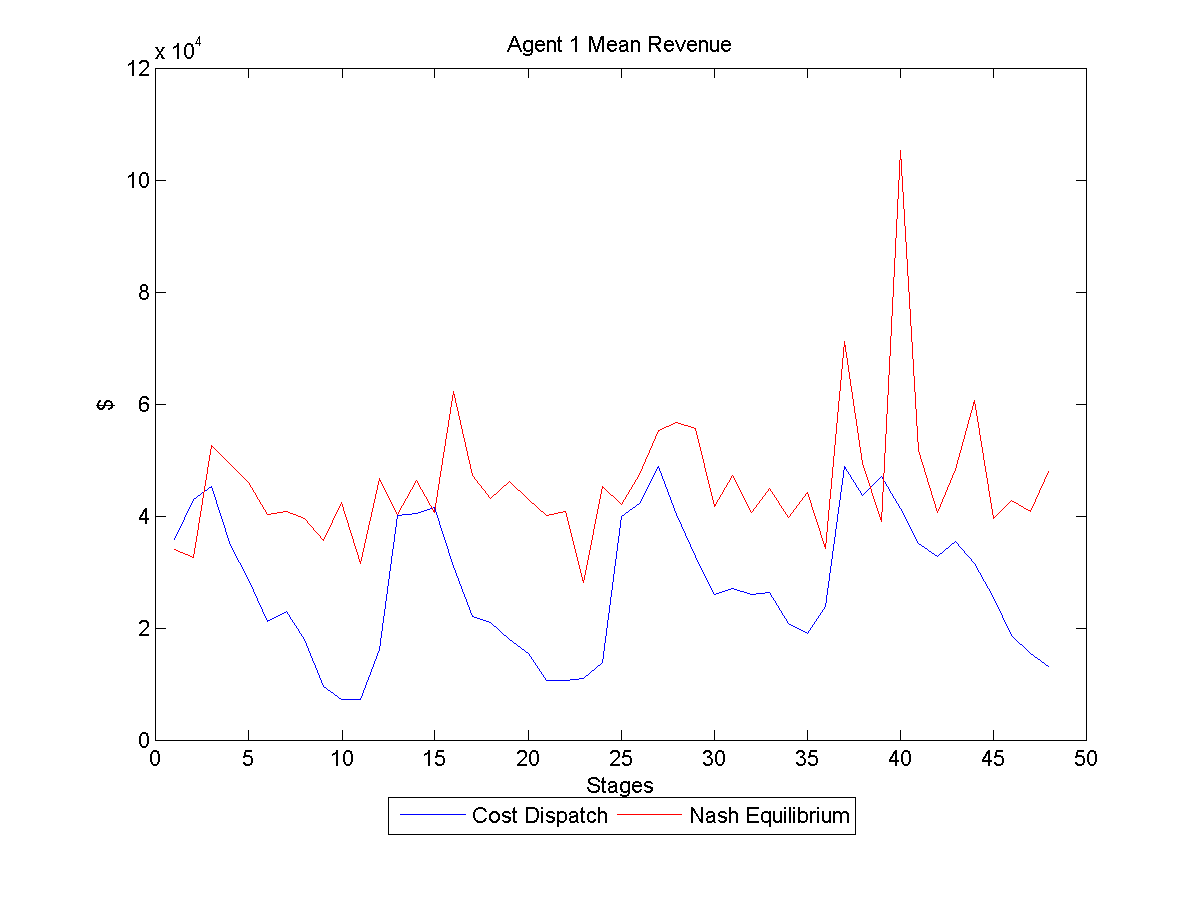}
\caption{Mean Revenue: Agent 1}\label{figAgentRev1}
\end{figure}

\begin{figure}[h!]
\centering
\includegraphics[width=0.7\linewidth]{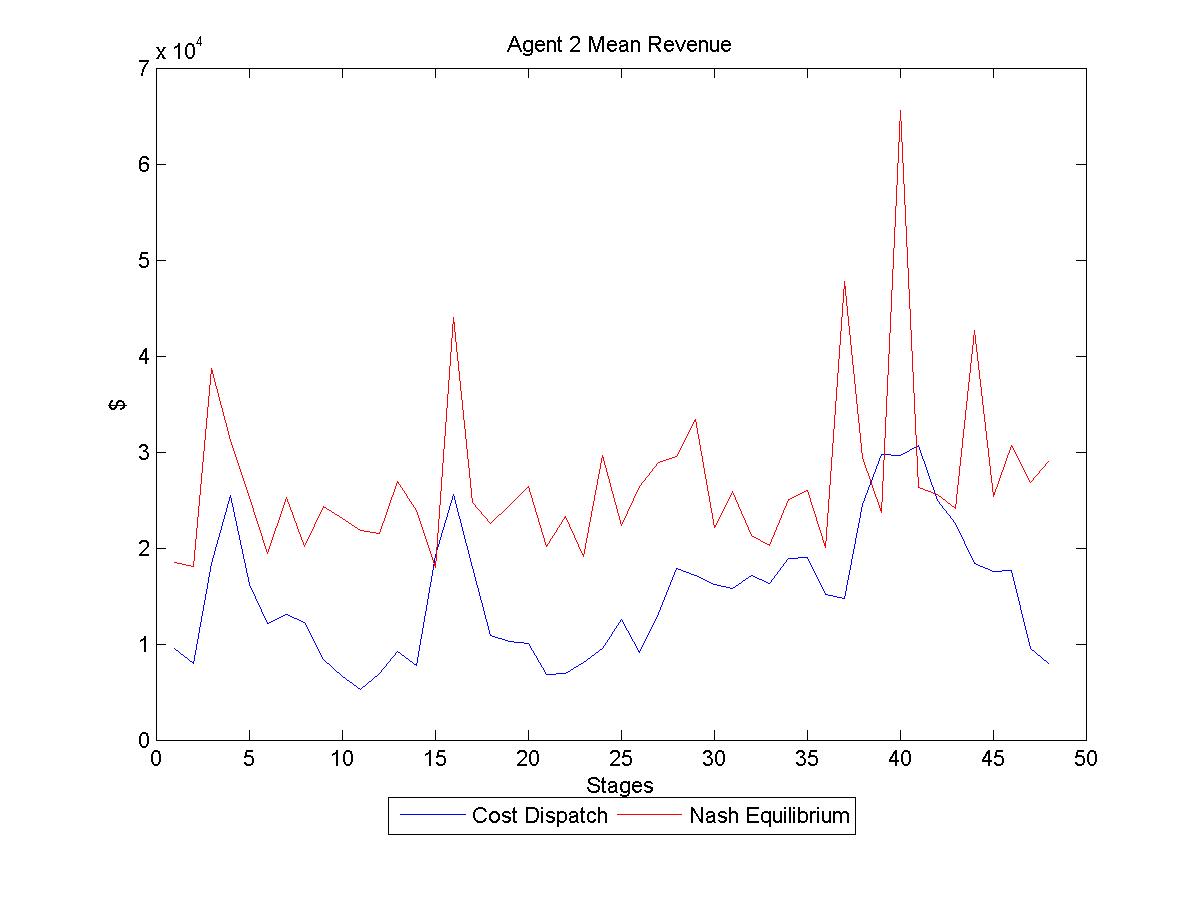}
\caption{Mean Revenue: Agent 2}\label{figAgentRev2}
\end{figure}

\begin{figure}[h!]
\centering
\includegraphics[width=0.7\linewidth]{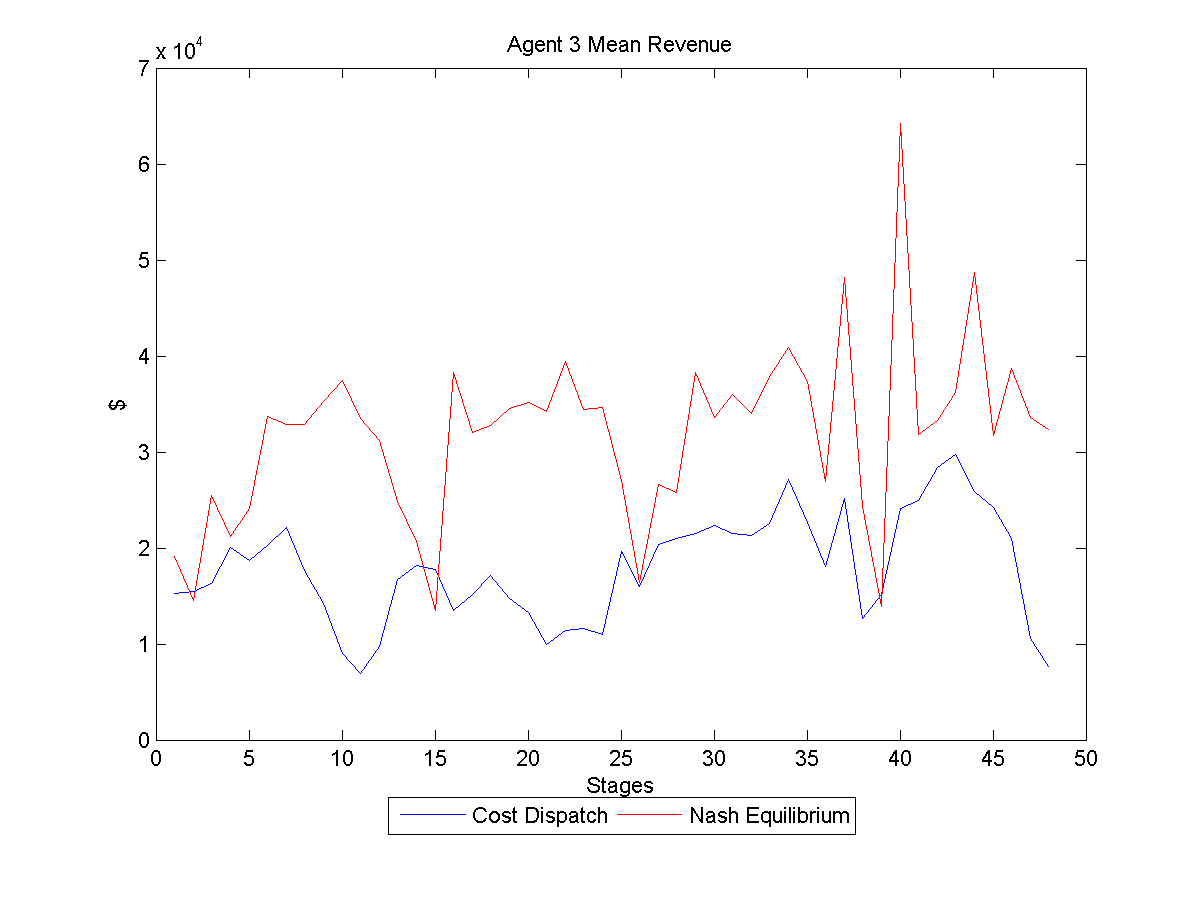}
\caption{Mean Revenue: Agent 3}\label{figAgentRev3}
\end{figure}

Fig. \ref{figAgentRev} displays the agents revenue resulting from both simulations: Centralized dispatch (CD) and Nash Equilibrium of deregulated market (NE). We can observe that all the agents indeed had a higher revenue in almost all stages. Fig. \ref{figAgentRev1}, Fig. \ref{figAgentRev2} and Fig. \ref{figAgentRev3} show the revenue result for each agent separately.

\clearpage

\begin{figure}[h!]
\centering
\includegraphics[width=0.8\linewidth]{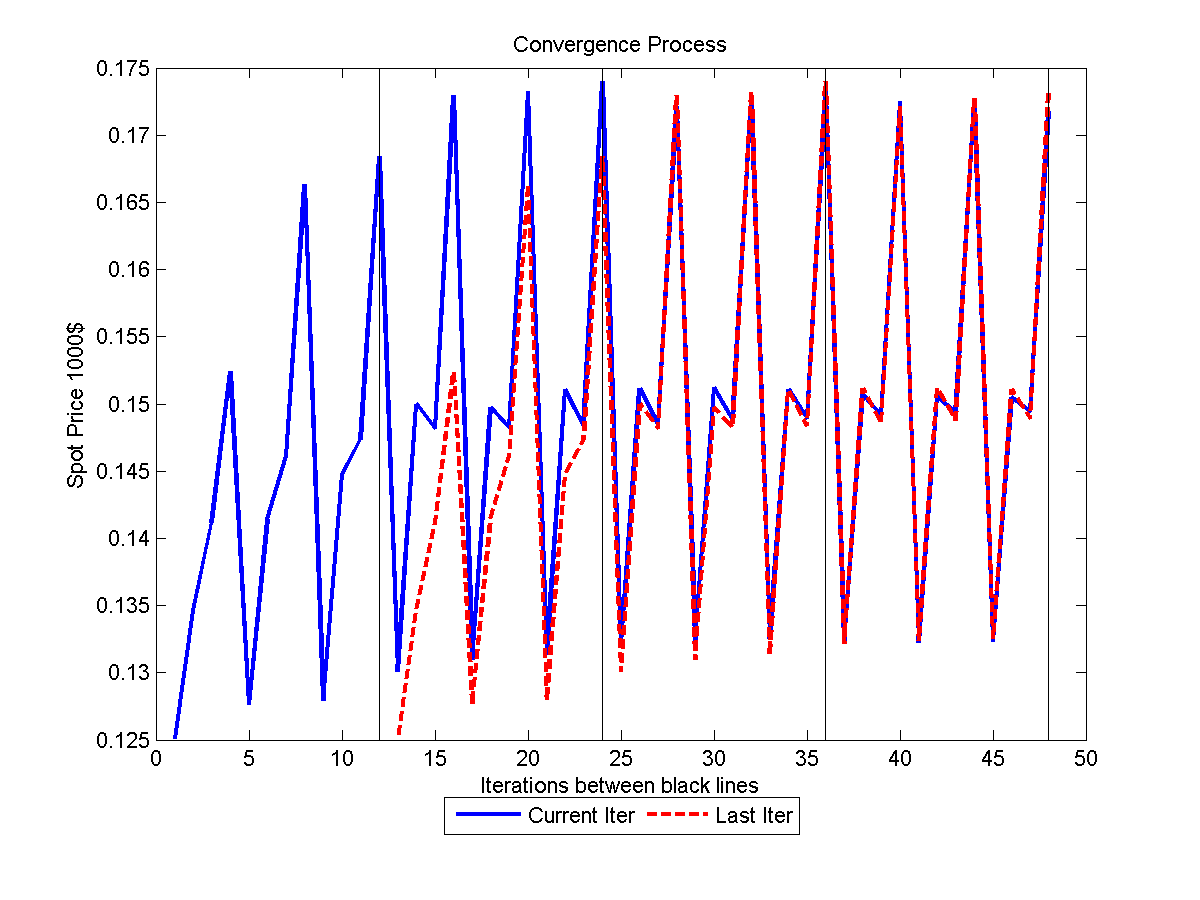}
\caption{Convergence in spot prices}\label{figSpotConv}
\end{figure}

Finally, we present Fig. \ref{figSpotConv} which depicts the convergence process in terms of the spot price average. The points in the graph are yearly average spot price. Each sequence of four points represents four consecutive yearly spot prices after one agents has optimized its policy, the following four points represent the yearly spot price after the second agent has optimized its policy. Thus, each sequence of twelve points represents a complete fixed point iteration through all the agents, for clarity we separated the iteration with vertical black lines. In the figure, the blue line presents the described that, while the red dashed line presents the same data shifted by twelve points, or one complete fixed point iteration, and we observe that ate the fourth complete iteration the process has stabilized.

\clearpage

\newpage

\section{Conclusions}

As we have been arguing since the very beginning of this work the problem of simulating power markets in presence of Hydro plants and Price Maker agents able to exercise market power is extremely difficult due to the complex combination of time-coupling and game theoretic equilibrium formulations.

We were successful in presenting a complete and detailed methodology to simulate a power market under those hydro-thermal and price maker assumptions. The simulation methodology made wide use of existing works and led to interesting results in our application to the Panama power system.

As expected we obtained a situation in which agents were able to modify the spot prices leading to more profitable situations than the centralized dispatch.


\subsection{Further Steps}

Some interesting questions of both theoretical an practical nature can be placed about this research area. Concerning game theory, we can ask in which situation does the game has multiple equilibria points or even what are the conditions for the existence of some equilibria.

The application of the methodology to other power systems with different number of price maker agents is another interesting research area. The questions of how to model agents which have hydro plants in the same cascade is extremely important and should be addressed in future works.

\clearpage

\bibliographystyle{ieeetr}
\bibliography{market_models}

\begin{thebibliography}{10}

\bibitem{hu2007using}
X.~Hu and D.~Ralph, ``Using epecs to model bilevel games in restructured
  electricity markets with locational prices,'' {\em Operations research},
  vol.~55, no.~5, pp.~809--827, 2007.

\bibitem{hunt2002making}
S.~Hunt, {\em Making competition work in electricity}, vol.~146.
\newblock John Wiley \& Sons, 2002.

\bibitem{al2006electricity}
A.~Al-Sunaidy and R.~Green, ``Electricity deregulation in oecd (organization
  for economic cooperation and development) countries,'' {\em Energy}, vol.~31,
  no.~6, pp.~769--787, 2006.

\bibitem{maurer2011electricity}
L.~Maurer and L.~A. Barroso, {\em Electricity auctions: an overview of
  efficient practices}.
\newblock World Bank Publications, 2011.

\bibitem{zhu1997review}
J.~Zhu and M.-y. Chow, ``A review of emerging techniques on generation
  expansion planning,'' {\em Power Systems, IEEE Transactions on}, vol.~12,
  no.~4, pp.~1722--1728, 1997.

\bibitem{barroso2006nash}
L.~A. Barroso, R.~D. Carneiro, S.~Granville, M.~V. Pereira, and M.~Fampa,
  ``Nash equilibrium in strategic bidding: a binary expansion approach,'' {\em
  Power Systems, IEEE Transactions on}, vol.~21, no.~2, p.~629, 2006.

\bibitem{pereira1991multi}
M.~V. Pereira and L.~M. Pinto, ``Multi-stage stochastic optimization applied to
  energy planning,'' {\em Mathematical Programming}, vol.~52, no.~1-3,
  pp.~359--375, 1991.

\bibitem{steeger2014optimal}
G.~Steeger, L.~A. Barroso, and S.~Rebennack, ``Optimal bidding strategies for
  hydro-electric producers: A literature survey,'' {\em Power Systems, IEEE
  Transactions on}, vol.~29, no.~4, pp.~1758--1766, 2014.

\bibitem{hobbs2000strategic}
B.~F. Hobbs, C.~B. Metzler, and J.-S. Pang, ``Strategic gaming analysis for
  electric power systems: An mpec approach,'' {\em Power Systems, IEEE
  Transactions on}, vol.~15, no.~2, pp.~638--645, 2000.

\bibitem{flach2010long}
B.~Flach, L.~Barroso, and M.~Pereira, ``Long-term optimal allocation of hydro
  generation for a price-maker company in a competitive market: latest
  developments and a stochastic dual dynamic programming approach,'' {\em IET
  generation, transmission \& distribution}, vol.~4, no.~2, pp.~299--314, 2010.

\bibitem{krishna2009auction}
V.~Krishna, {\em Auction theory}.
\newblock Academic press, 2009.

\bibitem{borenstein2000diagnosing}
S.~Borenstein, J.~Bushnell, and F.~Wolak, ``Diagnosing market power in
  california's restructured wholesale electricity market,'' tech. rep.,
  National Bureau of Economic Research, 2000.

\bibitem{wolak2014effective}
F.~Wolak, ``An effective regulator is needed for new zealand electricity
  industry,'' tech. rep., New Zealand Herald, 2014.

\bibitem{villar2003hydrothermal}
J.~Villar and H.~Rudnick, ``Hydrothermal market simulator using game theory:
  assessment of market power,'' {\em Power Systems, IEEE Transactions on},
  vol.~18, no.~1, pp.~91--98, 2003.

\bibitem{kelman2001market}
R.~Kelman, L.~A.~N. Barroso, and M.~V.~F. Pereira, ``Market power assessment
  and mitigation in hydrothermal systems,'' {\em Power Systems, IEEE
  Transactions on}, vol.~16, no.~3, pp.~354--359, 2001.

\bibitem{rebennack2010handbook}
S.~Rebennack, P.~M. Pardalos, M.~V. Pereira, and N.~A. Iliadis, {\em Handbook
  of power systems II}.
\newblock Springer, 2010.

\bibitem{taylor2015convex}
J.~A. Taylor, {\em Convex optimization of power systems}.
\newblock Cambridge University Press, 2015.

\bibitem{pereira1985stochastic}
M.~Pereira and L.~Pinto, ``Stochastic optimization of a multireservoir
  hydroelectric system: a decomposition approach,'' {\em Water resources
  research}, vol.~21, no.~6, pp.~779--792, 1985.

\bibitem{yakowitz1982dynamic}
S.~Yakowitz, ``Dynamic programming applications in water resources,'' {\em
  Water resources research}, vol.~18, no.~4, pp.~673--696, 1982.

\bibitem{kauppi2008empirical}
O.~Kauppi, M.~Liski, {\em et~al.}, {\em An empirical model of imperfect dynamic
  competition and application to hydroelectricity storage}.
\newblock Citeseer, 2008.

\bibitem{gjelsvik2010long}
A.~Gjelsvik, B.~Mo, and A.~Haugstad, ``Long-and medium-term operations planning
  and stochastic modelling in hydro-dominated power systems based on stochastic
  dual dynamic programming,'' in {\em Handbook of Power Systems I}, pp.~33--55,
  Springer, 2010.

\bibitem{gross2000generation}
G.~Gross and D.~Finlay, ``Generation supply bidding in perfectly competitive
  electricity markets,'' {\em Computational \& Mathematical Organization
  Theory}, vol.~6, no.~1, pp.~83--98, 2000.

\bibitem{bunn2000forecasting}
D.~W. Bunn, ``Forecasting loads and prices in competitive power markets: The
  technology of power system competition,'' {\em Proceedings of the IEEE},
  vol.~88, no.~2, pp.~163--169, 2000.

\bibitem{kwon2012optimization}
R.~H. Kwon and D.~Frances, ``Optimization-based bidding in day-ahead
  electricity auction markets: A review of models for power producers,'' in
  {\em Handbook of Networks in Power Systems I}, pp.~41--59, Springer, 2012.

\bibitem{conejo2001mathematical}
A.~J. Conejo and F.~J. Prieto, ``Mathematical programming and electricity
  markets,'' {\em Top}, vol.~9, no.~1, pp.~1--22, 2001.

\bibitem{ventosa2005electricity}
M.~Ventosa, A.~Baillo, A.~Ramos, and M.~Rivier, ``Electricity market modeling
  trends,'' {\em Energy policy}, vol.~33, no.~7, pp.~897--913, 2005.

\bibitem{boyd2004convex}
S.~Boyd and L.~Vandenberghe, {\em Convex optimization}.
\newblock Cambridge university press, 2004.

\bibitem{luo1996mathematical}
Z.-Q. Luo, J.-S. Pang, and D.~Ralph, {\em Mathematical programs with
  equilibrium constraints}.
\newblock Cambridge University Press, 1996.

\bibitem{pereira2005strategic}
M.~V. Pereira, S.~Granville, M.~H. Fampa, R.~Dix, and L.~A. Barroso,
  ``Strategic bidding under uncertainty: a binary expansion approach,'' {\em
  Power Systems, IEEE Transactions on}, vol.~20, no.~1, pp.~180--188, 2005.

\bibitem{hobbs2004modeling}
B.~F. Hobbs and U.~Helman, ``Modeling prices in competitive electricity
  markets: Chapter 3 complementarity-based equilibrium modeling for electric
  power markets,'' 2004.

\bibitem{conejo2002optimal}
A.~J. Conejo, J.~Contreras, J.~M. Arroyo, and S.~De~la Torre, ``Optimal
  response of an oligopolistic generating company to a competitive pool-based
  electric power market,'' {\em Power Systems, IEEE Transactions on}, vol.~17,
  no.~2, pp.~424--430, 2002.

\bibitem{de2003stackelberg}
M.~de~Luj{\'a}n~Latorre and S.~Granville, ``The stackelberg equilibrium applied
  to ac power systems-a noninterior point algorithm,'' {\em Power Systems, IEEE
  Transactions on}, vol.~18, no.~2, pp.~611--618, 2003.

\bibitem{cau2002co}
T.~D.~H. Cau and E.~J. Anderson, ``A co-evolutionary approach to modelling the
  behaviour of participants in competitive electricity markets,'' in {\em Power
  Engineering Society Summer Meeting, 2002 IEEE}, vol.~3, pp.~1534--1540, IEEE,
  2002.

\bibitem{baillo2004optimal}
A.~Baillo, M.~Ventosa, M.~Rivier, and A.~Ramos, ``Optimal offering strategies
  for generation companies operating in electricity spot markets,'' {\em Power
  Systems, IEEE Transactions on}, vol.~19, no.~2, pp.~745--753, 2004.

\bibitem{gjelsvik1999algorithm}
A.~Gjelsvik, M.~M. Belsnes, and A.~Haugstad, ``An algorithm for stochastic
  medium-term hydrothermal scheduling under spot price uncertainty,'' in {\em
  Proceedings of 13th Power Systems Computation Conference}, 1999.

\bibitem{fosso1999generation}
O.~B. Fosso, A.~Gjelsvik, A.~Haugstad, B.~Mo, and I.~Wangensteen, ``Generation
  scheduling in a deregulated system. the norwegian case,'' {\em Power Systems,
  IEEE Transactions on}, vol.~14, no.~1, pp.~75--81, 1999.

\bibitem{redondo1999short}
N.~J. Redondo and A.~Conejo, ``Short-term hydro-thermal coordination by
  lagrangian relaxation: solution of the dual problem,'' {\em Power Systems,
  IEEE Transactions on}, vol.~14, no.~1, pp.~89--95, 1999.

\bibitem{dieu2009improved}
V.~N. Dieu and W.~Ongsakul, ``Improved merit order and augmented lagrange
  hopfield network for short term hydrothermal scheduling,'' {\em Energy
  Conversion and Management}, vol.~50, no.~12, pp.~3015--3023, 2009.

\bibitem{catalao2012optimal}
J.~Catalao, H.~Pousinho, and J.~Contreras, ``Optimal hydro scheduling and
  offering strategies considering price uncertainty and risk management,'' {\em
  Energy}, vol.~37, no.~1, pp.~237--244, 2012.

\bibitem{pousinho2012scheduling}
H.~M. Pousinho, V.~M.~F. Mendes, and J.~P. Catal{\~a}o, ``Scheduling of a hydro
  producer considering head-dependency, price scenarios and risk-aversion,''
  {\em Energy Conversion and Management}, vol.~56, pp.~96--103, 2012.

\bibitem{mizyed1992operation}
N.~R. Mizyed, J.~C. Loftis, and D.~G. Fontane, ``Operation of large
  multireservoir systems using optimal-control theory,'' {\em Journal of Water
  Resources Planning and Management}, vol.~118, no.~4, pp.~371--387, 1992.

\bibitem{lohndorf2010optimal}
N.~L{\"o}hndorf and S.~Minner, ``Optimal day-ahead trading and storage of
  renewable energies an approximate dynamic programming approach,'' {\em Energy
  Systems}, vol.~1, no.~1, pp.~61--77, 2010.

\bibitem{helseth2010long}
A.~Helseth, B.~Mo, and G.~Warland, ``Long-term scheduling of hydro-thermal
  power systems using scenario fans,'' {\em Energy Systems}, vol.~1, no.~4,
  pp.~377--391, 2010.

\bibitem{lino2003bid}
P.~Lino, L.~A.~N. Barroso, M.~V. Pereira, R.~Kelman, and M.~H. Fampa,
  ``Bid-based dispatch of hydrothermal systems in competitive markets,'' {\em
  Annals of Operations Research}, vol.~120, no.~1-4, pp.~81--97, 2003.

\bibitem{scott1996modelling}
T.~J. Scott and E.~G. Read, ``Modelling hydro reservoir operation in a
  deregulated electricity market,'' {\em International Transactions in
  Operational Research}, vol.~3, no.~3-4, pp.~243--253, 1996.

\bibitem{baslis2011mid}
C.~G. Baslis and A.~G. Bakirtzis, ``Mid-term stochastic scheduling of a
  price-maker hydro producer with pumped storage,'' {\em Power Systems, IEEE
  Transactions on}, vol.~26, no.~4, pp.~1856--1865, 2011.

\bibitem{pousinho2013risk}
H.~M. Pousinho, J.~Contreras, A.~G. Bakirtzis, and J.~P. Catalao,
  ``Risk-constrained scheduling and offering strategies of a price-maker hydro
  producer under uncertainty,'' {\em Power Systems, IEEE Transactions on},
  vol.~28, no.~2, pp.~1879--1887, 2013.

\bibitem{wehinger2013modeling}
L.~A. Wehinger, G.~Hug-Glanzmann, M.~D. Galus, and G.~Andersson, ``Modeling
  electricity wholesale markets with model predictive and profit maximizing
  agents,'' {\em Power Systems, IEEE Transactions on}, vol.~28, no.~2,
  pp.~868--876, 2013.

\bibitem{thome2013non}
F.~Thome, M.~Pereira, S.~Granville, and M.~Fampa, ``Non-convexities
  representation on hydrothermal operation planning using sddp,'' tech. rep.,
  working paper, available: http://www. psr-inc. com.
  br/portal/psr/publicacoes, 2013.

\bibitem{philpott2016midas}
A.~Philpott, F.~Wahid, and B.~Fr{\'e}d{\'e}ric, ``Midas: A mixed integer
  dynamic approximation scheme,'' {\em Optimization-online}, 2016.

\bibitem{powell2007approximate}
W.~B. Powell, {\em Approximate Dynamic Programming: Solving the curses of
  dimensionality}, vol.~703.
\newblock John Wiley \& Sons, 2007.

\bibitem{sheble1994unit}
G.~B. Sheble and G.~N. Fahd, ``Unit commitment literature synopsis,'' {\em
  Power Systems, IEEE Transactions on}, vol.~9, no.~1, pp.~128--135, 1994.

\bibitem{street2014energy}
A.~Street, A.~Moreira, and J.~M. Arroyo, ``Energy and reserve scheduling under
  a joint generation and transmission security criterion: An adjustable robust
  optimization approach,'' {\em Power Systems, IEEE Transactions on}, vol.~29,
  no.~1, pp.~3--14, 2014.

\bibitem{dantzig1998linear}
G.~B. Dantzig, {\em Linear programming and extensions}.
\newblock Princeton university press, 1998.

\bibitem{bertsimas1997introduction}
D.~Bertsimas and J.~N. Tsitsiklis, {\em Introduction to linear optimization},
  vol.~6.
\newblock Athena Scientific Belmont, MA, 1997.

\bibitem{chvatal1983linear}
V.~Chvatal, {\em Linear programming}.
\newblock Macmillan, 1983.

\bibitem{bertsekas1995dynamic}
D.~P. Bertsekas, D.~P. Bertsekas, D.~P. Bertsekas, and D.~P. Bertsekas, {\em
  Dynamic programming and optimal control}, vol.~1.
\newblock Athena Scientific Belmont, MA, 1995.

\bibitem{bellman2015applied}
R.~E. Bellman and S.~E. Dreyfus, {\em Applied dynamic programming}.
\newblock Princeton university press, 2015.

\bibitem{birge2011introduction}
J.~R. Birge and F.~Louveaux, {\em Introduction to stochastic programming}.
\newblock Springer Science \& Business Media, 2011.

\bibitem{hamilton1994time}
J.~D. Hamilton, {\em Time series analysis}, vol.~2.
\newblock Princeton university press Princeton, 1994.

\bibitem{durbin2012time}
J.~Durbin and S.~J. Koopman, {\em Time series analysis by state space methods}.
\newblock No.~38, Oxford University Press, 2012.

\bibitem{mcleod1994diagnostic}
A.~I. McLeod, ``Diagnostic checking of periodic autoregression models with
  application,'' {\em Journal of Time Series Analysis}, vol.~15, no.~2,
  pp.~221--233, 1994.

\bibitem{bezerra2012assessment}
B.~Bezerra, A.~Veiga, L.~A. Barroso, and M.~Pereira, ``Assessment of parameter
  uncertainty in autoregressive streamflow models for stochastic long-term
  hydrothermal scheduling,'' in {\em Power and Energy Society General Meeting,
  2012 IEEE}, pp.~1--8, IEEE, 2012.

\bibitem{puterman2014markov}
M.~L. Puterman, {\em Markov decision processes: discrete stochastic dynamic
  programming}.
\newblock John Wiley \& Sons, 2014.

\bibitem{bertsekas1976dynamic}
D.~P. Bertsekas, ``Dynamic programming and stochastic control,'' {\em
  Mathematics in science and engineering}, vol.~125, pp.~222--293, 1976.

\bibitem{ross2014introduction}
S.~M. Ross, {\em Introduction to probability models}.
\newblock Academic press, 2014.

\bibitem{benders1962partitioning}
J.~F. Benders, ``Partitioning procedures for solving mixed-variables
  programming problems,'' {\em Numerische mathematik}, vol.~4, no.~1,
  pp.~238--252, 1962.

\bibitem{shapiro2011analysis}
A.~Shapiro, ``Analysis of stochastic dual dynamic programming method,'' {\em
  European Journal of Operational Research}, vol.~209, no.~1, pp.~63--72, 2011.

\bibitem{philpott2008convergence}
A.~B. Philpott and Z.~Guan, ``On the convergence of stochastic dual dynamic
  programming and related methods,'' {\em Operations Research Letters},
  vol.~36, no.~4, pp.~450--455, 2008.

\bibitem{homem2011sampling}
T.~Homem-de Mello, V.~L. de~Matos, and E.~C. Finardi, ``Sampling strategies and
  stopping criteria for stochastic dual dynamic programming: a case study in
  long-term hydrothermal scheduling,'' {\em Energy Systems}, vol.~2, no.~1,
  pp.~1--31, 2011.

\bibitem{de2015improving}
V.~L. De~Matos, A.~B. Philpott, and E.~C. Finardi, ``Improving the performance
  of stochastic dual dynamic programming,'' {\em Journal of Computational and
  Applied Mathematics}, vol.~290, pp.~196--208, 2015.

\bibitem{mas1995microeconomic}
A.~Mas-Colell, M.~D. Whinston, J.~R. Green, {\em et~al.}, {\em Microeconomic
  theory}, vol.~1.
\newblock Oxford university press New York, 1995.

\bibitem{nash1951non}
J.~Nash, ``Non-cooperative games,'' {\em Annals of mathematics}, pp.~286--295,
  1951.

\bibitem{fudenberg1991game}
D.~Fudenberg and J.~Tirole, ``Game theory, 1991,'' {\em Cambridge,
  Massachusetts}, vol.~393, 1991.

\bibitem{osborne1994course}
M.~J. Osborne and A.~Rubinstein, {\em A course in game theory}.
\newblock MIT press, 1994.

\bibitem{nash1950equilibrium}
J.~F. Nash {\em et~al.}, ``Equilibrium points in n-person games,'' {\em Proc.
  Nat. Acad. Sci. USA}, vol.~36, no.~1, pp.~48--49, 1950.

\bibitem{glicksberg1952further}
I.~L. Glicksberg, ``A further generalization of the kakutani fixed point
  theorem, with application to nash equilibrium points,'' {\em Proceedings of
  the American Mathematical Society}, vol.~3, no.~1, pp.~170--174, 1952.

\bibitem{cournot1838recherches}
A.-A. Cournot, {\em Recherches sur les principes math{\'e}matiques de la
  th{\'e}orie des richesses par Augustin Cournot}.
\newblock chez L. Hachette, 1838.

\bibitem{bertrand1883theorie}
J.~Bertrand, ``Th{\'e}orie math{\'e}matiques de la richesse sociale,'' {\em
  Journal des Savants}, vol.~48, pp.~499--508, 1883.

\bibitem{nisan2007algorithmic}
N.~Nisan, T.~Roughgarden, E.~Tardos, and V.~V. Vazirani, {\em Algorithmic game
  theory}, vol.~1.
\newblock Cambridge University Press Cambridge, 2007.

\bibitem{von2007theory}
J.~Von~Neumann and O.~Morgenstern, {\em Theory of games and economic behavior}.
\newblock Princeton university press, 2007.

\end{thebibliography}

\clearpage

\section{Appendix 1 - Hydro Thermal Power Systems Operation}\label{hydrothermal}

	The state of art solution of the hydro thermal dispatch problem is given by the SDDP algorithm, which first appeared in the seminal paper by Pereira \textit{et al.} \cite{pereira1991multi}, many other methods and variants were proposed to solve the operation problem such as the approximate dynamic programming \cite{powell2007approximate}. However we stick to SDDP because it is still on the leading frameworks to solve the problem and its heavily used in industrial applications.
	
	Aiming to understand this method and build the cornerstone of the following algorithm the SDDP will be carefully described. The SDDP will be built over a sequence of simpler steps and classical methods of mathematical programming. Firstly the thermal operation will be described in Section \ref{thermalop}, in the sequence, Section \ref{htop} aims to point out the main difficulties that arise from considering hydro plants in the system operation. Before the advent of SDDP the problem of Section \ref{htop} would be solved by Stochastic Dynamic Programming which has its own merits and is easily applied to a much wider class of problems, in chapter \ref{market} we present a combination of SDDP and SDP used to solve an specific problem. Finally the SDDP method is presented with some modifications from the one first shown in \cite{pereira1991multi}.

\subsection{Thermal Operation}\label{thermalop}

	The operation of thermal system has a lot of nuances of its own and its still a extremely important study field. A lot of effort has been put to solve many specific thermal operation problems. Some of these very interesting problem arise from the combination of an extremely simple problem with the most varied characteristics of real life power systems such as the Unit-Commitment problem \cite{sheble1994unit}, its contigency constrained counterpart \cite{street2014energy}, the optimal power flow \cite{taylor2015convex}, thermal expansion \cite{}, network expansion \cite{}.
	
	The simplest of the thermal operation problems considers that a system us fully described by its demand $d$ and a set $G$ of thermal plants. Moreover, each thermal plant $i\in G$ is described by its operation cost $c_i$ and its generation capacity $\overline{g_i}$. Such representation of a power system yield to the following thermal dispatch problem on the variable $g_i$:
	
\begin{align}
& {\text{minimize}} && \sum_{i \in G}{c_j {g}_{i}} \\
& \text{subject to} \quad \notag \\
&&& \sum_{i \in G}{g_{i}}  = d_{}  && \leftarrow {\pi_p}_{} && \label{simpleloadbal}\\
&&& g_{i} \leq \overline{g}_i && &&, \forall i \in G 
\end{align}
	
	The solution of this problem is the minimum cost to meet the given demand and the values of $g_i$ in the optimum, represented as $g_i^*$. Will be the optimal dispatch. From now on all the variables with superscript $*$ represent their values at the optimal solution. The above problem is a standard optimisation problem that belongs to the class of linear programming \cite{dantzig1998linear}, these problems can be solved by classical methods such as the Simplex Method and Barrier Methods \cite{bertsimas1997introduction}\cite{chvatal1983linear}. 
	
	In the above problem, ${\pi_p}_{}$ represents the dual variable that corresponds to the first constraint (\ref{simpleloadbal}). Dual variables are fundamental concepts of optimization theory \cite{boyd2004convex}, also known by economists as shadow prices since they represent the infinitesimal variation of the optimal problem value as the right hand side of some constraint is varied.
	
	Its easily demonstrated that the thermal problem of this section can be solved by a straightforward algorithm in which all the plants are ordered and we choose the cheaper ones to attain the demand, the dual variable will be the cost of the most expensive plant dispatch. We call this dispatch the merit ordered dispatch.
	
	To finish this section we highlight two facts. Firstly, we could have multi-period thermal dispatch which is fundamental in unit-commitment for instance. However under our simplification hypothesis which only take into account cost, capacity and demand the problem has no temporal coupling and therefore the multi-period operation decouple into single period operation. Secondly, this problem is the archetype of a much broader class of problems which include any kind of market clearing where agents are represented by pairs of price and quantity.
	
\subsection{Hydro-Thermal Operation}\label{htop}

	Including hydro plants into power system operation brings to it two main difficulties: time coupling and stochastic operation. These two characteristics make the problem extremely harder to solve when compared to the simple thermal dispatch of last section. In this section we describe these two difficulties so that they can be dealt with in the next sections algorithms.
	
	Hydro plants typically include water reservoirs, which are used as storage mechanisms that can be used to save incoming water instead of using it as it comes. The exception are the so called run of the river hydros which do not have reservoirs, thus they turbine water as it comes. This capacity of saving water  to forthcoming periods creates a temporal coupling in the hydro-thermal operation problem. 
	
	This coupling can be used to reduce the operation cost of the power system. Hydro plant have no operation cost since the \textit{fuel} used in power generation is water, differently from thermal plants which must buy some fuel such as Oil, Diesel, Carbon etc. Indeed both plants have some other costs such as maintenance, but the classical simplification does not consider these costs. Therefore one could conclude that since hydro plants have cost zero they should always come before thermals in dispatch order. This is wrong.
	
	In order to understand how saving water is important consider the following example. Imagine a system composed of two thermal plants and one hydro,that will be operated in 2 stages. The thermals have capacity 10MW and 15MW and cost 50\$ by MW and 200\$ by MW. The demand is always 20MW. The hydro has a full reservoir that enables it to produce 20MW, and no water will arrive to the reservoir. Since the hydro have zero cost one could use all the water in the first stage and have a zero cost operation, but then, in the second stage the cost of the system would be $10 \times 50 + 10 \times 200 = 2500$. However, if we decided to use only half of the water in the first stage so we could use the other half on the other stage our total operation cost would be $2\times 10 \times 50 = 1000$. Therefore saving water can be extremely important in the operation problem.
	
	Mathematical problems with period coupling are frequently solved using  a well known method called Dynamic Programming \cite{bertsekas1995dynamic} which is commonly used in control and graph problems. In the core of this method is the construction of a recursive equation known as the Bellman equation \cite{bellman2015applied}. We will define a Bellman-equation for the hydro-thermal operation problem in the next sections.
	
	Considering hydro plants will take the problem away from the deterministic situation to the stochastic realm. The problem is that one can never know how much water will arrive in some reservoir because it depends on physical phenomena such as the rain and ice defrost.
	
	The incoming water of the reservoirs actually has two different sources, water used of upstream power plants and natural water, the second one will be called inflow from now on. The inflows are represented in the problem by a discrete time stochastic processes that materializes in random variables for each hydro plant and stage. This takes our problem to the field of Stochastic Programming	 \cite{birge2011introduction}. 
	
	Its extremely hard to obtain analytical solutions to general  stochastic programs, one of the common methods to solve stochastic programs is the so called \textit{scenario approach}. In this method, many realizations of the underlying stochastic processes are sampled from a priori fixed distributions as in a Monte Carlo method. Finally the problem is solved for all those scenarios jointly to reach an approximate solution of the problem.
	
	Inflows are stochastic processes that typically have some underlying temporal structure. This processes are known as time series and form a whole study field of its own theoretical and practical importance \cite{hamilton1994time}. Many possible models can be used to capture the temporal and spatial structure of the inflow process the most famous are the AR processes \cite{hamilton1994time}, more recently state-space models have also been used to estimate linear time series models \cite{durbin2012time}. The canonical model used in power systems operation is the Periodic Autoregressive model known as PAR \cite{mcleod1994diagnostic}, which will be extremely useful in the next sections because is retains convexity  of the model and captures periodic structure inherent to inflows process. Many Studies have devoted attention to the use and choice of such models \cite{bezerra2012assessment}.
	
	The scenario approach to solve stochastic problems usually takes use of the concept of a scenario tree \cite{} which is a sample of the stochastic process. The solution of such problem is the optimal operation of the system. However, the computational burden to solve the problem using simply scenario trees and extend formulations grows very rapidly with the problem size. In order to solve this exponentially growing problem Dynamic Programming is applied as a decomposition method.

	\subsection{Stochastic Dynamic Programming}

Literature in stochastic dynamic programming, SDP, is vast and some general bookshelf references include \cite{puterman2014markov}, \cite{bertsekas1976dynamic} and \cite{ross2014introduction}. Before the success of \cite{pereira1991multi} SDP was the fundamental technique to solve the dispatch, but it still leads to combinatorial explosion.

In the SDP approach the multi-stage and multi-scenario problem is decomposed in single stage deterministic (sub-)problems, that is the problem for each stage and (inflow) scenario is solved separately. The complete solution is given by solving the problems in reverse chronological order: we start by solving the stage $t$ sub-problems (one for each scenario), the solution information is passed to stage $t-1$ by means of the aforementioned Bellman recursion. In other words the problem for some stage is solved in a finite set of scenarios and from these values we construct a function to represent the cost at stage $t$ in stage $t-1$, hence we name it a \textit{future cost function}. Note that we can presents arbitrary functions by simply using piecewise linear interpolations and MILP for instance.

The time coupling in  such problem is due to the inflow process and the reservoir condition, therefore the information linking to stages is the cost of the system operation as a function of reservoir level and conditioned inflow values. These time coupling variables are known as \textit{State Variable}.

A standard hydro thermal dispatch problem, with reservoir and inflow states, has the following recursion:

\begin{align}
&\alpha_t^{SDP}(\hat{v}_{t,H}^s,[\hat{a}_{t,H}^l]) =\notag\\ 
& {\text{minimize}} \ \ \ \sum_{i \in G, b \in B}{c_j {g}_{t,b,i}} + \frac{1}{|L|}\sum_{l \in L}{\alpha_{t+1}^{SDP}({v}_{t+1,H},[{a}_{t+1,H}^l])} \label{sdpObj}\\
& \text{subject to} \quad & & \notag \\
& v_{t+1,i} = \hat{v}_{t,i}^s + \hat{a}_{t,i}^s - (u_{t,i} + x_{t,i}) + \sum_{j \in M(i)}(u_{t,j} + x_{t,j}) & & \leftarrow {\pi_h}_{t,i}& &,  \forall i \in H  \label{sdpHydroBal}   \\
& \sum_{b \in B}{e_{t,b,i}} = \rho_i u_{t,i} && &&, \forall i \in H \label{sdpHydroEne} \\ 
& \sum_{i \in H}{e_{t,b,i}} + \sum_{i \in G}{g_{t,b,i}}  = d_{t,b} - \sum_{i \in R}{r_{t,b,i}}  && \leftarrow {\pi_p}_{t,b} &&, \forall b \in B \label{sdpLoadBal}\\
& v_{t+1,i} \leq \overline{v}_i && &&, \forall i \in H \\
& u_{t,i} \leq \overline{u}_i && &&, \forall i \in H \\
& e_{t,b,i} \leq \overline{e}_i && &&, \forall i \in H \\ 
& g_{t,b,i} \leq \overline{g}_i && &&, \forall i \in G \\
& a_{t+1,i}^l = \sum_{p \in P(i,t)} {\phi_{t,i}^p \hat{a}_{t+1-p,i}^s }+ \hat{\xi}_{t,i}^l  && &&, \forall i \in H, l \in L\label{sdpAR}
\end{align}

Firstly we have to highlight that this is the sub-problem from a single stage and scenario, therefore the subscript $t$ is fixed as well as the superscript $s$. In this optimization problem we have as decision variables: the generation of thermal plant $i$ in block $b$ given by $g_{t,b,i}$, the spillage, turbined water, and reservoir level at the end of the stage from hydro $i$ given, respectively by $u_{t+1,i}$, $x_{t,i}$ and $v_{t,i}$; and the hydro energy produced in each block by each plant $e_{t,b,i}$. 

Considering blocks is efficient approximation for the so called load levels for instance. In typical power systems, demand varies significantly with hours of the day and days of weeks, the load in weekends or at night is way smaller.

The sets in the above defined mathematical program are: $H$ the set of Hydro plants, $G$ the set of thermal plants, $R$ the set of renewable energy plants, $B$ the set of blocks, $M(i)$ the set of hydros that spill and turbine to the reservoir of hydro $i$, finally we ha the set $L$ of openings. The set of openings is just an extra way of accounting for variability, so given the inflow scenario $t$ we can obtain $|L|$ possible outcomes for stage $t+1$, instead of just one, and represent future variability. This sophistication is extremely useful for modelling Markov-chain SDDP models.

Finally we have the following constraints in the sub-problem: (\ref{sdpHydroBal}
) represents the hydro balance; (\ref{sdpHydroEne}) enforces that the energy produced by each hydro in all blocks must match the turbined water; (\ref{sdpLoadBal}) represents the load balance that enforces how much energy must be produced to meet the demand and finally (\ref{sdpAR}) describes the autoregressive process of inflows. 

The objective function (\ref{sdpObj}) is composed of two terms: the \textit{Immediate Cost Function} given by the thermal production costs and the \textit{Future Cost Function} that represents the water values since its a function of inflows and storage in stage $t+1$. The function $\alpha_{t+1}^{SDP}({v}_{t+1,H},[{a}_{t+1,H}^l])$ is constructed from the solution of the problems in stage $t+1$ and can be represented by MILP as we previously commented.

We also have two important dual variables that come from the solution of some sub-problem: ${\pi_h}_{t,i}$ which is the water marginal cost for plant $i$ and the system's energy spot price ${\pi_p}_{t,b}$ in block $b$.

In a standard SDP algorithm we initialize by fixing the states for which we shall solve problems to gather points and construnct the function $$\alpha_{t+1}^{SDP}({v}_{t+1,H},[{a}_{t+1,H}^l])$$. Therefore we construct a discrete set $M$ of storage vectors, $\hat{v}_{t,H}= \hat{v}_{t,H}^1, \dots, \hat{v}_{t,H}^M$ and another set $L$ of inflow vectors $\hat{a}_{t,H}= \hat{a}_{t,H}^1, \dots, \hat{a}_{t,H}^L$. This procedure is done for all stages and one must note that for the inflows the inflows are conditioned on the ones from the  previous stages, due to the process temporal structure. 
	
Now we have the following algorithm for the SDP version of this problem:
	
\begin{algorithm}[H]
	\caption{SDP for Hydro Thermal Operation}
	\label{Algorithm} 
	\begin{algorithmic}
		\STATE{initialize discrete sets $\hat{v}_{t,H}$ and $\hat{a}_{t,H}, \forall t \in T$ }
		\FOR{$t = |T|, |T|-1, \dots, 2, 1$}
		\FORALL{$\hat{v}_{t,H}^i \in \{ \hat{v}_{t,H}^1, \dots, \hat{v}_{t,H}^M \}$}
		\FORALL{$\hat{a}_{t,H}^i \in \{ \hat{a}_{t,H}^1, \dots, \hat{a}_{t,H}^L \}$}
		\STATE{solve problem (\ref{sdpObj})-(\ref{sdpAR})}
		
		\ENDFOR
		\ENDFOR
		\STATE{Create future cost function for stage $t-1$ : $\alpha_{t}^{SDP}({v}_{t,H},[{a}_{t,H}^l])$}
		\ENDFOR
	\end{algorithmic}
\end{algorithm}	

In this algorithm the number of problems that must be solve is of the order $|T|\times|M|\times|L|$. The main problem of such method is that a system may have multiple hydros and a reasonable approximation of their reservoir levels and inflows would require a great combination of possible values. If we fix $K$ possible states for both inflows and storage values we have $M = L = K^{|H|}$, therefore the problem also explodes computationally even for a small number $K$.

This motivates the necessity for a next algorithm to reduce even more the computational burden of the solution of the Hydro Thermal operation.

\subsection{Stochastic Dynamic Dual Programming}

Stochastic Dual Dynamic Programming as presented in \cite{pereira1991multi}  makes clever use of duality theory to approximate the Bellman equation without having to create huge discrete sets and iterate through them.
	
	SDP approximates a function by evaluating it in a discrete and finite number of points, in case the function value is needed in some different point some interpolation must be carried on. For convex problems its easier, while for non-convex its typically hard and one must rely on Mixed Integer Linear Programming, MILP.
	
	SDDP rely on the convexity of the problem to improve the approximation of the Bellman equation known as \textit{Future Cost Function} (or \textit{Future Benefit Function} in the case of maximization problems) in the standard SDDP nomenclature. The fundamental idea is to use, not only information of value of previously evaluated functions, but also its derivative. In other words, while SDP performs a zero order approximation of an arbitrary function by evaluating it at some points, SDDP performs a first order approximation of the \textit{convex} future cost function by evaluating a function and its derivative at some points. The derivative can be obtained by relying on dual multipliers in a way extremely similar to the Benders Decomposition \cite{benders1962partitioning}.
	
	The future cost function therefore must be convex in the standard SDDP since its epigraph shall be approximated by a set of hyperplanes that, in theory, can represent a convex function with arbitrary precision. Those hyperplanes are constructed only from the function evaluation and its derivative. The epigraph of the future cost function can be described by a set of hyperplanes and the functions itself is given by the following optimization problem:
	
\begin{align}
& \alpha_{t+1}({v}_{t,H},[{a}_{t+1,H}^l]) = \hspace{0.1in}\notag \\
& {\text{minimize}} \ \ \ \ \alpha_{t+1} \label{fcfOBJ_2}\\
& \text{subject to} \quad \alpha_{t+1} \geq \varphi^{m}_{t+1} + \sum_{i \in H}{{\varphi_h}_{t+1,i}^m  {v}_{t,i}}+ \sum_{i \in H, p \in P(i,t)}{{\varphi_a}_{t+1,1,i}^m {a}_{t+1-p+1,i}^l} &&, \forall m \in \mathscr{M}, l \in L \label{fcfHyperPlanes_2}
\end{align}

Given such approximation, the subproblem for each stage and scenarios and SDDP will be basically the same as the SDP one, by simply replacing the future cost function:

\begin{align}
&\alpha_t(\hat{v}_{t,H}^s,[\hat{a}_{t,H}^s]) = \notag\\
& {\text{minimize}} & & \sum_{i \in G, b \in B}{c_j {g}_{t,b,i}} + \frac{1}{|L|}\sum_{l \in L}{\alpha_{t+1}({v}_{t+1,H},[{a}_{t+1,H}^l])} \label{sddpOBJ_2} \\
& \text{subject to} \quad & & v_{t+1,i} = \hat{v}_{t,i}^s + \hat{a}_{t,i}^s - (u_{t,i} + x_{t,i}) + \sum_{j \in M(i)}(u_{t,j} + x_{t,j}) & & \leftarrow {\pi_h}_{t,i}& &,  \forall i \in H    \\
&&& \sum_{b \in B}{e_{t,b,i}} = \rho_i u_{t,i} && &&, \forall i \in H \\ 
&&& \sum_{i \in H}{e_{t,b,i}} + \sum_{i \in G}{g_{t,b,i}}  = d_{t,b} - \sum_{i \in R}{r_{t,b,i}}  && \leftarrow {\pi_p}_{t,b} &&, \forall b \in B \\
&&& v_{t+1,i} \leq \overline{v}_i && &&, \forall i \in H \\
&&& u_{t,i} \leq \overline{u}_i && &&, \forall i \in H \\
&&& e_{t,b,i} \leq \overline{e}_i && &&, \forall i \in H \\ 
&&& g_{t,b,i} \leq \overline{g}_i && &&, \forall i \in G \\
&&& a^l_{t+1,i} = \sum_{p \in P(i,t)} {\phi_{t,i}^p \hat{a}_{t+1-p,i}^s }+ \hat{\xi}_{t,i}^l  && &&, \forall i \in H, l \in L \label{sddpAR_2}
\end{align}

	By simply adding the equations (\ref{fcfHyperPlanes_2}) to the problem (\ref{sddpOBJ_2})-(\ref{sddpAR_2}) we have a stage and scenario sub-problem completely represented as a linear program. 
	
	Using derivatives imply in a much better approximation, however we will not have a significant computational gain if points must be evaluated in a discrete pre-determined set as in SDP, at this point the second big idea comes into play. Once more following the steps of the Benders decomposition, the future cost function will be approximated iteratively by a scheme similar to the master/slave approach. In other words, instead of generating an arbitrary number of hyperplanes we iteratively generate only the \textit{interesting} ones.
	
	The algorithm will be briefly explained and then detailed in the forthcoming paragraphs. We start by defining initial approximations of the future cost functions in all stages by empty sets of hyperplanes. Secondly, one define a number of scenarios and define initial values for state variables, typically the current situation of reservoirs and inflows. The innovation vectors $\hat{\xi}_{t,i}^l$, sampled from the pertinent distributions and known as random shocks or residuals, will define the variability between inflow scenarios in the autoregressive process and, consequently, generate some scenario tree. These innovation vectors can be the determined a priori or on-the-run, typically we have a set $L$ of shocks whose cardinality is smaller than the number o scenarios. The number $|L|$ is called the number of openings, since each shock correspond to a tree opening, the set $S$ of scenarios are commonly known as \textit{forward series}.

	Having initialized the algorithm and defined the scenario update rule, the inflow AR equation(\ref{sddpAR_2}), we be begin to discuss the core algorithm. One defines feasible storage values for all stages and scenarios then starting in the last stage $t = |T|$  the operation is optimized for all scenarios in $S$, in the last stage there is no future cost function as in the SDP algorithm. The solution of each scenario gives a hyperplane to approximate the future cost function of stage $t-1$, one proceeds this way until solving all scenarios up to stage $t = 2$, this part of the algorithm is known as \textit{Backward Recursion}.
	
	In the second part of the algorithm, know as \textit{Forward Simulation}, the operation problem is solved in chronological order. Starting from stage $t = 1$, solutions of each stage problem (together with the autoregressive update rule) yield to updates in the values of storage and inflow state variables so that a locally optimal operation is obtained. From these new values of state variables another Backward Recursion takes place generating new hyperplanes and improving the future cost function approximation. We proceed this way until convergence is met.
	
	The convergence of the algorithm is said to be achieved if the average total cost of the first stage objective value stays in some confidence interval of the average sum of all stages \textit{Immediate cost function}. The Immediate Cost function is simply the total objective function of some stage minus its future cost function.
	
	The SDDP algorithm is summarized in the next lines:
	
\begin{algorithm}[H]
	\caption{SDDP for Hydro Thermal Operation}
	\label{Algorithm2} 
	\begin{algorithmic}
	\WHILE{Converge not met}
		\STATE{Forward Simulation:}
		\FOR{$t = 1, \dots, |T|-1$}
		\FORALL{$s \in S$}
		\STATE{Solve problem (\ref{sddpOBJ_2})-(\ref{sddpAR_2})}
		\STATE{Update the reservoir levels $v_{t+1,i}$ of stage $t+1$}
	 	\ENDFOR
		\ENDFOR
		\STATE{Backward Recursion:}
		\FOR{$t = 1, \dots, |T|$}
		\FORALL{$s \in S$}
		\STATE{solve problem (\ref{sddpOBJ_2})-(\ref{sddpAR_2})}
		\STATE{Generate cut for stage $t-1$ from objective value, solution value ad dual variables}
	 	\ENDFOR
		\ENDFOR
	\ENDWHILE
	\end{algorithmic}
\end{algorithm}

Because of its great practical success the SDDP algorithm has been studied continuously. The method has been deeply analysed in light of stochastic programming framework \cite{shapiro2011analysis}, the method convergence have been analysed by \cite{shapiro2011analysis} and \cite{philpott2008convergence}, stopping criteria have been also studied in \cite{homem2011sampling}, performance has been studied in \cite{de2015improving}.

\section{Appendix 2 - Game Theory Overview}

Game theory is a common framework for studying economic systems where not all agents can be assumed price takers and therefore market power must be modelled. It can be seen as an extension of the simpler social welfare maximization theory\cite{mas1995microeconomic}. 

More specifically we shall rely on the basics of non-cooperative game theory\cite{nash1951non} where agents aim to optimize their own utility functions, \cite{fudenberg1991game}\cite{osborne1994course}, which is basically a measure of happiness, without accounting for coalitions. In simple words, agents optimize their utility functions by choosing appropriate strategies. In this work utilities will be simply revenues and strategies will be bids, either quantity bids or price bids.

We say that we have a \textit{Nash equilibrium} when, given a set of bids(strategies) for all agents, no single agents benefits(improve their revenues) from deviating from such strategy. Basics and variants of Nash equilibrium including existence theorems have been studied since the seminal work of John Nash \cite{nash1950equilibrium}. Many existence theorems and method for attaining them are based on fixed-point theorems\cite{glicksberg1952further}.

Two more concepts that shall be useful in the sequence are the ones of Cournot and Bertrand models. The Cournot model\cite{cournot1838recherches} assumes the agents compete in the market via their production level, in our case the energy quantity offer. On the other hand, in the Bertrand model\cite{bertrand1883theorie} agents compete via their prices, which for us are the bid prices. These two models represent, respectively, the least and most efficient markets when taken to the limit case of a duopoly where Cournot prices are some function of the demand and Bertrand quantities are a single value per agent, where the best price takes all the demand.

Of course game theory has its own limitations the two more important were commented previously on the text: the algorithmic difficulties which is a fertile research area \cite{nisan2007algorithmic} and the rationality hypothesis that come from the expected utility hypothesis \cite{von2007theory}, which raises the question of whether real life organizations clever enough to be players in a game.

%
%


\end{document}